\documentclass{amsart}
\usepackage{amstext,amssymb,amsbsy,array,amsmath,txfonts,color,todonotes}

\renewcommand{\div}{\text{div$\,$}}
\newcommand{\Div}{\text{{\bf div}$\,$}}

\newcommand{\bfv}{{\boldsymbol v}}

\newcommand{\bfsig}{{\boldsymbol\sigma}}

\newcommand{\bfp}{{\boldsymbol p}}
\newcommand{\bfP}{{\boldsymbol P}}
\newcommand{\bfy}{{\boldsymbol y}}

\newcommand{\bfI}{{\boldsymbol I}}
\newcommand{\bfu}{{\boldsymbol u}}
\newcommand{\bbR}{{\mathbb{R}}}
\newcommand{\IR}{{\mathbb{R}}}

\newcommand{\triangulation}{\mathcal{T}}
\newcommand{\mcF}{\mathcal{F}}
\newcommand{\mcT}{\mathcal{T}}
\newcommand{\mcP}{\mathcal{P}}

\newcommand{\mcS}{\mathcal{S}}

\newcommand{\bfn}{{\boldsymbol n}}

\newcommand{\bfeps}{{\boldsymbol\varepsilon}}

\newcommand{\bft}{{\boldsymbol t}}

\newcommand{\bfx}{{\boldsymbol x}}

\newcommand{\tn}{|\mspace{-1mu}|\mspace{-1mu}|}

\newcommand{\Eb}{E_\Omega}
\newcommand{\nub}{\nu_\Omega}
\newcommand{\Es}{E_\Sigma}
\newcommand{\tb}{t_\Omega}

\newcommand{\betab}{\beta_\Omega}
\newcommand{\nablas}{\nabla_\Sigma}

\newcommand{\mcCP}{\mathcal{C}_P}
\newcommand{\mcCB}{\mathcal{C}_B}

\newtheorem{rmk}{Remark}

\begin{document}
\title{A Simple Approach for Finite Element Simulation of Reinforced Plates}
\author[EB]{Erik Burman}
\author[PH]{Peter Hansbo}
\author[ML]{Mats G. Larson}
\address[EB]{Department of Mathematics, University College London, London, UK--WC1E 6BT, United Kingdom} 
\address[PH]{Department of Mechanical Engineering, J\"onk\"oping University, S-551 11 J\"onk\"oping,
Sweden}
\address[ML]{Department of Mathematics and Mathematical Statistics,
Ume{\aa} University, SE--901 87 Ume{\aa}, Sweden}
\date{}
\begin{abstract}
We present a new approach for adding Bernoulli beam reinforcements to Kirchhoff plates. 
The plate is discretised using a continuous/discontinuous finite element method based on 
standard continuous piecewise polynomial finite element spaces. The beams are discretised 
by the CutFEM technique of letting the basis functions of the plate represent also the 
beams which are allowed to pass through the plate elements. This allows for a fast and easy 
way of assessing where the plate should be supported, for instance, in an optimization loop.
\end{abstract}
\maketitle

\section{Introduction}

Reinforcements of solids using lower--dimensional structures such as beams can be simulated in a finite element  context by
coupling the variables of the beam to the variables of the solid, either along element edges as in McCune, Armstrong, and Robinson \cite{McArRo00} or by interpolation on element edges as in Sadek and Shahrour \cite{SaSh04}. In the latter case, the beam geometry can be modelled independently of the bulk mesh which is crucial; however, the finite element approximation of the 
lower--dimensional object is otherwise independent and uncoupled to the solid, and the rotation degrees of freedom of beams are hard to match to the solid (if they are to influence the solution in the solid).

\begin{figure}
\begin{center}
\includegraphics[scale=0.3]{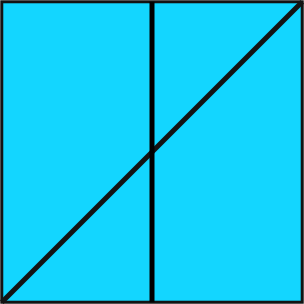}
\hspace{1cm}
\includegraphics[scale=0.3]{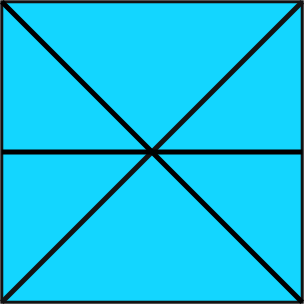}
\hspace{1cm}
\includegraphics[scale=0.3]{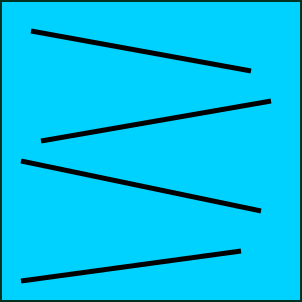}
\end{center}
\caption{Examples of plates reinforced by beams.}
\end{figure}

\textcolor{black}{In \cite{HaLaLa17} we proposed to use the same finite element space for the beam as for the higher dimensional structure; more precisely, the trial and test space for the beam is obtained by taking the restriction or trace to the beam. Here we further develop this approach to allow for coupling between plates and beams, more precisely  the Kirchhoff-Love plate model and the Bernoulli beam model.} These models involve 
fourth order partial differential equations. We discretize these models using the so called continuous/discontinuous Galerkin, c/dG, method which relaxes the required $C^1$ continuity of the
 shape functions for the beam and plate by use of a discontinuous Galerkin approach with $C^0$--continuity. \textcolor{black}{We emphasise that the concept is quite general, as illustrated in our previous work 
on embedding in elastic solids, of membranes \cite{CeHaLa16} and of embedded trusses and beams \cite{HaLaLa17}.} A 
similar approach was recently suggested for modelling embedded trusses by L\'e, Legrain, and Mo\"es \cite{LeLeMo17}.

\section{Modeling of Reinforced Plates}

\subsection{The Basic Approach}
In this Section we develop a simple model of a set of beam elements in a plate. The 
main approach is as follows: 
\begin{itemize}
\item Given a continuous finite element space, based on at least 
second order polynomials for the plate, we define the finite element space for the 
one--dimensional structure as the restriction of the plate finite element 
space to the structure which is geometrically modeled by an embedded 
curve or line.

\item To formulate a finite element method on the restricted or trace finite 
element space we employ continuous/discontinuous Galerkin approximations 
of the Euler--Bernoulli beam model. 
The beams are then modeled using the CutFEM paradigm and the 
stiffness of the embedded beams is in the most basic version, which we 
consider here, simply added to the plate stiffness. 

\end{itemize}

To ensure coercivity of the cut beam model we in general need to add a certain stabilization 
term which provides control of the discrete functions variation in the vicinity of the beam. 
However, for beams embedded in a plate, the plate stabilizes the beam discretizations, and 
we shall show that if the plate is stiff enough compared to the beam the usual additional 
stabilization \cite{BurClaHan15} is superfluous. The plate problem may also be viewed as 
an interface problem in order to more accurately approximate the plate in the vicinity of the 
beam structure; this approach is however significantly more demanding from an implementation 
point of view and we leave it for future work. 

The work presented here is an extension of earlier work \cite{CeHaLa16}
where membrane structures were considered, in which case a linear approximation 
in the bulk suffices.

\begin{figure}
\begin{center}
\includegraphics[scale=0.4]{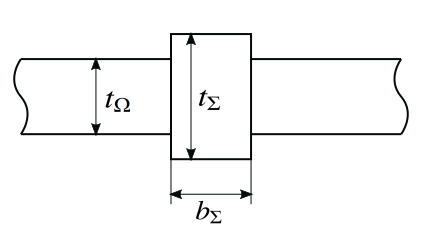}
\hspace{0.5cm}
\includegraphics[scale=0.4]{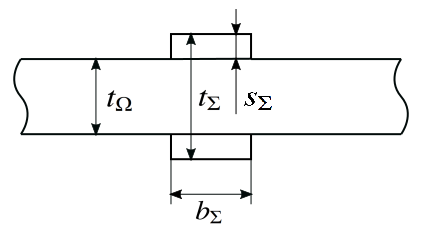}
\end{center}
\caption{Left: The reinforced plate geometry parameters, $t_\Omega$, $t_\Sigma$, 
and $b_\Sigma$. Right: Alternative design of reinforcement with two separate beams 
of thickness $s_\Sigma = (t_\Sigma - t_\Omega)/2$ above and below the plate.}\label{fig:plate-geometry}
\end{figure}

\subsection{The Kirchhoff--Love Plate Model}

In the Kirchhoff--Love plate model, posed on a polygonal domain 
$\Omega\subset \bbR^2$ with boundary $\partial\Omega$ and exterior unit normal 
$\bfn$, we seek an out--of--plane 
(scalar) displacement $u$ to which we associate the strain (curvature) tensor
\begin{equation}
\bfeps(\nabla u) := \frac12\left(\nabla\otimes (\nabla u) + (\nabla u ) \otimes \nabla \right) 
= \nabla \otimes \nabla u = \nabla^2 u   
\end{equation}
and the plate stress (moment) tensor
\begin{align}\label{eq:plate-stress-tensor}
\bfsig_P (\nabla u) &:= \mcCP \left(\bfeps(\nabla u) + \nub (1- {\nub })^{-1}
\div\nabla u \, \bfI \right)
\\
&= \mcCP \left( \nabla^2 u + \nub (1-\nub)^{-1} \Delta u \bfI \right)
\end{align}
where 
\begin{equation}\label{eq:mcCP}
\mcCP =  \frac{\Eb \tb^3}{12(1+\nub)} 
\end{equation}
with $\Eb$ the Young's modulus, $\nub$ the Poisson's ratio, and $\tb$ denotes the 
plate thickness. Since $0 \leq  \nub \leq 0.5$ the constants are uniformly bounded. 

The Kirch\-hoff--Love problem then takes the form: given the out--of--plane 
load (per unit area) $f$, find the displacement $u$ such that
\begin{align}
\div \Div \bfsig_P ( \nabla u )= f &  \qquad \text{in $\Omega$}
\\
u = 0  & \qquad \text{on $\partial \Omega$}
\\
\bfn \cdot \nabla u = 0   &\qquad \text{on $\partial \Omega$}
\end{align}
where $\Div$ and $\text{div}$ denote the divergence of a tensor and a vector 
field, respectively.

\paragraph{Weak Form}
The variational problem takes the form: 
Find the displacement $u \in V_\Omega = H^2_0(\Omega)$ such that
\begin{equation}\label{eq:varform}
a_\Omega(u, v )= l_\Omega (v) \qquad \forall v \in V_\Omega
\end{equation}
where the forms are defined by 
\begin{align}
a_\Omega(v,w) & = (\bfsig_P(\nabla v), 
\bfeps(\nabla w))_\Omega 
\\
l_\Omega (v ) &= (f,v)_\Omega
\end{align}

We employ the following notation: 
$L^2(\omega)$ is the Lebesgue space of square integrable 
functions on $\omega$ with scalar product 
$(\cdot,\cdot)_{L^2(\omega)} = (\cdot,\cdot)_{\omega}$
and  $(\cdot,\cdot)_{L^2(\Omega)} = (\cdot,\cdot)$, 
and norm  $\|\cdot \|_{L^2(\omega)} = \|\cdot\|_\omega$ 
and $\|\cdot\|_{L^2(\Omega)} = \|\cdot \|$,
$H^s(\omega)$ is the Sobolev space of order $s$ on 
$\omega$ with norm $\|\cdot\|_{H^s(\omega)}$,  
and $H^1_0(\Omega)=\{v\in \text{$H^1(\Omega):$}\; v  = 0 
\text{ on $\partial \Omega$}\}$, and $H^2_0(\Omega)=\{v\in H^2(\Omega):\, v = \bfn \cdot \nabla v = 0 
\text{ on $\partial \Omega$}\}$.

\subsection{The Euler--Bernoulli Beam Model}

Consider a straight thin beam with centerline $\Sigma \subset \Omega$ and 
a rectangular cross-section with width $b_\Sigma$ and thickness $t_\Sigma$, 
see Figure \ref{fig:plate-geometry}.
The modeling of the beam is performed using tangential differential calculus and 
we follow the exposition in \cite{HaLaLa14,HaLaLa17}, which also covers 
curved beams. Using this approach the beam equation is expressed in the 
same coordinate system as the plate, which is convenient in the construction 
of the cut finite element method for reinforced plates. 

Let $\bft$ be the tangent vector to the line $\Sigma$, embedded in $\IR^2$. 
We let $\bfp:\mathbb{R}^2\rightarrow \Sigma$ be the closest point mapping, 
i.e. $\bfp (\bfx ) = \bfy$ where $\bfy \in \Sigma$ minimizes the Euclidean norm 
$|  \bfx - \bfy |_{\IR^3}$. We define $\zeta$ as the signed distance function 
$\zeta(\bfx) := \pm\vert\bfx - \bfp(\bfx)\vert$, positive on one side of $\Sigma$ 
and negative on the other. 

Let $\bfP_\Sigma = \bft \otimes \bft$ be the projection onto the one dimensional 
tangent space of $\Sigma$ and define the tangential derivatives 
\begin{equation}\label{eq:tangential-derivatives}
\nablas v = \bfP_\Sigma \nabla v,\qquad \partial_{\bft} v = \bft \cdot \nabla v
\end{equation}
Then we have the identity 
\begin{equation}\label{eq:tangential-der-relation}
\nablas v = (\partial_{\bft} v) \bft
\end{equation}

Based on the assumption that planar cross sections orthogonal 
to the midline remain plane
after deformation we assume that the displacement takes the form 
\begin{equation} \label{displacementfield}
\bfu = u\bfn + \theta  \zeta \bft  
\end{equation}
where $\theta:\Sigma \rightarrow \mathbb{R}$ is an angle representing 
an infinitesimal rotation, assumed 
constant in the normal plane. In Euler--Bernoulli beam theory the beam 
cross-section is assumed plane and orthogonal to the beam midline after 
deformation and no shear deformations occur, which means that we have
\begin{align}
\theta = \bft\cdot\nabla u := \partial_{\bft}u
\end{align}
This definition for $\theta$ in combination with \eqref{displacementfield} 
constitutes the Euler--Bernoulli kinematic assumption 
\begin{align*}
\bfu = u \bfn + \zeta (\partial_{\bft} u )\bft = u \bfn + \zeta \nablas u
\end{align*}

We assume the usual Hooke's law for one dimensional structural members
\begin{align} \label{constitutive}
\bfsig_\Sigma(\bfu) = E_\Sigma \bfeps_\Sigma(\bfu)
\end{align}
where $\Es$ is the Young modulus and the tangential strain 
tensor is given by
\begin{equation}
\bfeps_\Sigma(\bfu) 
=\bfP_\Sigma\bfeps(\bfu) \bfP_\Sigma 
= \zeta \bfeps_\Sigma (\nablas u  ) 
%
\end{equation}
where in the last equality we used the identity 
\begin{equation}
\bfu \otimes \nabla 
=
(u\bfn + \zeta \nablas u ) \otimes \nabla
=
 \bfn \otimes (\nabla u) 
+ \zeta (\nablas u )\otimes \nabla 
\end{equation}
to conclude that 
\begin{equation}
\bfeps_\Sigma(\bfu) = \zeta \bfeps_\Sigma(\nablas u )
\end{equation}
%
%

Next note that the strain energy density can  be written
\begin{equation}
\bfsig_\Sigma(\bfu) : \bfeps_\Sigma(\bfu) 
= 
\zeta^2\bfsig_\Sigma(\nablas u ): \bfeps(\nablas u) 
\end{equation}
and the total energy of the beam structure is obtained by 
integrating over the beam volume
\begin{equation}
{\mathcal E}_\Sigma= \frac12\int_\Sigma
 I_\Sigma   \bfsig_\Sigma (\nablas u ) : \bfeps (\nablas u )    \, d\Sigma
 -\int_\Sigma a_\Sigma  f_\Sigma u\, d\Sigma
\end{equation}
where the integral over the cross section is accounted for by the cross-section area 
and its second moment
\begin{equation}
a_\Sigma = b_\Sigma t_\Sigma, \qquad I_\Sigma =b_\Sigma t_\Sigma^3/12
\end{equation}

We are thus led to introducing the beam stress tensor
\begin{equation}\label{eq:beam-stress-tensor}
\bfsig_{B,\Sigma}(\nablas v ) 
= I_\Sigma \bfsig_\Sigma (\nablas v)
=  I_\Sigma E_\Sigma \bfeps_\Sigma (\nablas v)
\end{equation}
and thus we have the beam Hooke law
\begin{equation}
\bfsig_{B,\Sigma}(\nablas v ) = \mcCB \bfeps_\Sigma (\nablas v)
\end{equation}
where 
\begin{equation}\label{eq:mcCB}
\mcCB = E_\Sigma I_\Sigma = \frac{E_\Sigma b_\Sigma t_\Sigma^3}{12} 
\end{equation}

Taking variations we obtain the weak statement, assuming zero displacements and 
rotations at the end points of $\Sigma$, we thus seek 
$u \in  V_\Sigma = H^2_0(\Sigma)$, such that
\begin{equation}
a_\Sigma(u,v) = l_\Sigma(v)\qquad \forall v \in V_\Sigma
\end{equation}
where the forms are defined by
\begin{equation}
a_\Sigma(v,w) =  \int_\Sigma
 \bfsig_{B,\Sigma} (\nablas v ) : \bfeps_\Sigma (\nablas w )  \, d\Sigma,
\qquad
l_\Sigma(v) = \int_\Sigma a_\Sigma  f_\Sigma  v  \, d\Sigma 
\end{equation}

\begin{rmk} We have the  identity  
\begin{align}\label{eq:identity-eps}
\bfeps_\Sigma (\nablas v ) = (\partial_{\bft}^2 v) \bft \otimes \bft 
\end{align}
since 
$(\nablas v) \otimes \nablas 
=  
((\partial_{\bft} v) \bft)\otimes \nablas 
= 
(\partial_{\bft} (\partial_{\bft} v) \bft)\otimes \bft
= 
(\partial^2_{\bft} v ) \bft \otimes \bft
$,
and thus
\begin{equation}\label{eq:identity-sig-eps}
\bfsig_{B,\Sigma}(\nablas v ):\bfeps(\nablas w) = \Es I_\Sigma \partial_{\bft}^2 v \partial_{\bft}^2 w
\end{equation}
which leads to 
\begin{equation}
a_\Sigma(v,w) 
= 
\int_\Sigma \bfsig_{B,\Sigma} (\nablas v ):\bfeps(\nablas w)\, d \Sigma
=
\int_\Sigma  \Es  I_\Sigma \partial_{\bft}^2 v \partial_{\bft}^2 w \, d \Sigma
\end{equation}
Here we recognize the right hand side as the traditional bilinear form associated with 
the Euler-Bernoulli beam. 
\end{rmk}
\textcolor{black}{
\begin{rmk} We note that in the alternative reinforcement geometry, right in Figure \ref{fig:plate-geometry}, we 
have 
\begin{equation}\label{eq:alternative-parameters}
a_\Sigma = b_\Sigma(t_\Sigma - t_\Omega),\qquad I_\Sigma = \frac{E_\Sigma b_\Sigma ( t^3_\Sigma - t^3_\Omega )}{12}
\end{equation}
We may also consider more complicated cross sections and compute the proper parameters.
\end{rmk}
}

\subsection{The Reinforced Plate Model}

Let $\mcS = \{ S \}$ be a set of beams arbitrarily oriented in $\Omega$. 
Using superposition we obtain the problem: find $u \in V$ such that 
\begin{equation}
a(u, v) = l(v)\qquad \forall v \in V
\end{equation}
where 
\begin{equation}
V = V_\Omega \bigcap_{\Sigma \in \mcS} V_\Sigma 
\end{equation}
and the forms are defined by
\begin{align}
a(v,w) & = a_\Omega(v,w) + \sum_{\Sigma \in \mcS} a_\Sigma(v,w)
\\
l(v) &= l_\Omega (v) + \sum_{\Sigma \in \mcS} l_\Sigma(v)
\end{align}

\textcolor{black}{
\begin{rmk}
Note that for the alternative plate reinforcement geometry, right in Figure \ref{fig:plate-geometry}, 
there is no geometric error in our method if we use the parameters 
(\ref{eq:alternative-parameters}).  In the 
standard reinforcement geometry, left in Figure \ref{fig:plate-geometry}, there is a 
however a geometric error proportional to $b_\Sigma$ in the plate bilinear form, which 
arises in the superposition since the intersection between the beam and the plate 
is nonempty. We will later see that $b_\Sigma$ typically is smaller 
(in practice significantly smaller) than the mesh size since we are using thin beam 
and plate theory, see (\ref{eq:Cmesh}), and thus the geometric error is small.
\end{rmk}
}

\section{Finite Element Discretization}

\subsection{The Mesh and Finite Element Spaces}

\begin{itemize}
\item We consider a subdivision $ \mcT_h =\{ T\}$ of $\Omega$ into
a geometrically conforming finite element mesh, with mesh parameter 
$h \in (0,h_0]$. We assume that the
elements are shape regular, i.e., the quotient of the diameter of the
smallest circumscribed sphere and the largest inscribed sphere is 
uniformly bounded. We denote by $h_T$ the diameter of element $T$ and by
$h = \max_{T \in \mcT_h} h_T$ the global mesh size parameter. 

\item Since we are using thin plate and beam theory we assume that 
there is a constant $C_{\text{mesh}}$ such that
\begin{equation}\label{eq:Cmesh}
C_{\text{mesh}}  \max(t_\Omega,t_\Sigma, b_\Sigma) \leq h
\end{equation} 

\item We shall use continuous, piecewise polynomial approximations, for 
both the membrane and plate problem. Let 
\begin{equation}
V_{\Omega,h,k} =\{ v \in C^0(\Omega):~\text{$v \vert_T \in \mcP_{k}(T)$
$\forall T \in \triangulation$}\}
\end{equation}
where $\mcP_k(T)$ is the space of polynomials of degree less or equal 
to $k$ defined on $T$. For simplicity, we write $V_{\Omega,h} = V_{\Omega,h,k}$.

\item To define our method we introduce the set of faces (edges) $F$ 
in the mesh, $\mcF_h =\{ F \}$, and we split $\mcF_h$ into two disjoint 
subsets
\begin{equation}
  \mcF_h = \mcF_{h,I} \cup \mcF_{h,B}
\end{equation}
where $\mcF_{h,I}$ is the set of faces in the interior of $\Omega$ and
$\mcF_{h,B}$ is the set of faces on the boundary. 

\item Further, with each
face $F$ we associate a fixed unit normal $\bfn_F$  such that for faces
on the boundary $\bfn_F$ is the exterior unit normal. We denote the
jump of a function $\bfv$ at a face $F$ by $ \left[\bfv
\right] = \bfv^+ - \bfv^-$ for $F \in \mcF_{h,I}$ and $
\left[\bfv\right] = \bfv^+$ for $F \in \mcF_{h,B}$, and the average
$\langle  \bfv\rangle = (\bfv^+ + \bfv^-)/2$ for $F \in \mcF_{h,I}$
and $\langle \bfv\rangle = \bfv^+$ for $F \in \mcF_{h,B}$, where
$\bfv^{\pm} = \lim_{\epsilon\downarrow 0} \bfv(\bfx\mp
\epsilon\,\bfn_F)$ with $\bfx\in F$.

\item Given a line segment $\Sigma$ in $\Omega$ that represents a beam 
we let $$\mcT_h(\Sigma) = \{ T \in \mcT_h : T \cap \Sigma \neq \emptyset \}$$
and we let $\mcF_h(\Sigma)$ be the set of all interior faces in $\mcT_h(\Sigma)$. 

\item The intersection points between $\Sigma$ and element faces in 
$\mcF_h(\Sigma)$ is denoted
\begin{equation}
\mcP_h (\Sigma) = \{ \bfx: \, \bfx = F \cap\Sigma, \,F \in \mcF_h(\Sigma) \} 
\end{equation}
and we assume that this is a discrete set of points (thus excluding 
the case where any $F\in\mcF_h$ coincides with a part of $\Sigma$). 

\begin{figure}
\begin{center}
\includegraphics[scale=0.3]{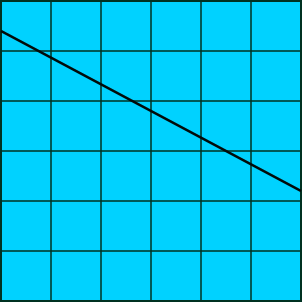}
\hspace{1cm}
\includegraphics[scale=0.3]{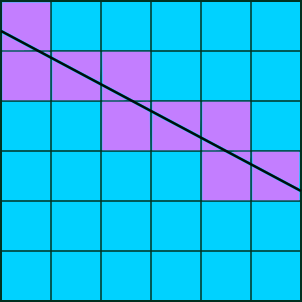}
\hspace{1cm}
\includegraphics[scale=0.3]{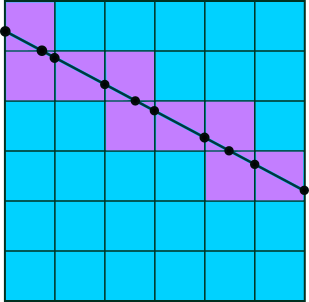}
\end{center}
\caption{The mesh $\mcT_h$ with one beam, the active mesh $\mcT_h(\Sigma)$ for the 
beam in purple, and the set of intersection points $\mcP_h(\Sigma)$.}
\end{figure}

\end{itemize}

\subsection{The c/dG Method for the Plate}

We approximate the solution to the plate problem using the continuous/discontinuous Galerkin 
(c/dG) method: Find $u_h \in V_{\Omega,h}$, with $k\geq 2$, such that
\begin{equation}\label{eq:cdg}
a_{\Omega,h}(u_h, v )   = l_\Omega(v) \qquad\forall v \in V_{\Omega,h} 
\end{equation}
The bilinear form $a_{\Omega,h}(\cdot,\cdot)$ is defined by
\begin{align} \label{nitsche_form}
a_{\Omega,h}(v,w) 
= {}& \sum_{T \in \mcT_h} ( \bfsig_P(\nabla v),\bfeps( \nabla w))_T
\\ \nonumber
&\qquad  - \sum_{F \in \mcF_{h,I} \cup \mcF_{h,B} }
 ( \langle\bfn_F \cdot \bfsig_P(\nabla v) \rangle,  [ \nabla w
  ])_F 
\\ \nonumber
&\qquad -  \sum_{F \in \mcF_{h,I} \cup \mcF_{h,B} }( [\nabla v ], \langle \bfn_F \cdot \bfsig_P(\nabla w)\rangle)_F
\\ \nonumber
&\qquad +  \sum_{F \in \mcF_{h,I} \cup
\mcF_{h,B}} \betab h_F^{-1} ([ \nabla v ], \left[ \nabla w\right])_F
\end{align}
\textcolor{black}{Here $\betab$ is a positive parameter of the form 
\begin{equation}
\betab =\beta_{\Omega,0} \mcCP = \beta_{\Omega,0} \frac{  \Eb \tb^3}{12(1+\nub)}
\end{equation}
where $\beta_{\Omega,0}$ is a constant depending on the polynomial order $k$,
see \cite{HanLar09} for details},  and $h_F$ 
is defined on each face $F$ by
\begin{equation}
h_F = \left( |T^+| + |T^-| \right) / ( 2 \, |F| ) \quad \text{for $F
= \partial T^+ \cap
\partial T^-$}
\end{equation}
with $|T|$ the area of $T$ and $| F|$ the length of $F$. 

\begin{rmk}
The idea of using 
continuous/discontinuous approximations was first proposed by Engel et al. \cite{EngGar02} 
and later analysed for Kirchhoff--Love and Mindlin--Reissner plates in 
\cite{HaLa03,HanLar09,HaHeLa11}, cf. also Wells and Dung \cite{WeDU07}.
\end{rmk}


\begin{rmk} Other boundary conditions for plates, for instance simply supported and free, 
can easily be included in the c/dG finite element method, see \cite{HanLar02} for details.  
\end{rmk}

\begin{rmk}\label{rem:bending-moments} For $v\in V_{\Omega,h}$ we have $[\nabla v] = [\bfn_F \cdot \nabla v] \bfn_F$ 
since $v$ is continuous across a face and $v=0$ on $\partial \Omega$. Therefore 
\begin{equation}
 ( \langle\bfn_F \cdot \bfsig(\nabla v) \rangle,  [ \nabla w
  ])_F  
  = ( \langle\bfn_F \cdot \bfsig(\nabla v) \cdot \bfn_F  \rangle, 
  [ \bfn_F \cdot \nabla w ])_F \quad 
\end{equation}
for all $v,w \in V_{\Omega,h}$, and we note that $(\bfn_F \cdot \bfsig(\nabla u)\cdot \bfn_F)|_F$ 
is the bending moment at the edge $F$.
\end{rmk}

\subsection{The Cut c/dG Method for a Beam}

We propose the following cut c/dG method. 
Find $u_h \in V_{\Sigma,h} = V_{\Omega,h} |_{\mcT_h(\Sigma)}$ such that
\begin{equation}
A_{\Sigma,h}(u_h,v) = l_{\Sigma}(v)\qquad \forall v \in  V_{\Sigma,h}
\end{equation}
where 
\begin{align}
A_{\Sigma,h}(v,w) & = a_{\Sigma,h}(v,w) + s_{\Sigma,h}(v,w)
\\
a_{\Sigma,h}(v,w) 
&=
\sum_{T\in \mcT_h(\Sigma)}  (\bfsig_{B,\Sigma}(\nablas v ),\bfeps(\nablas w ) )_{\Sigma\cap T} 
\\ \nonumber
&\qquad 
-\sum_{\bfx\in\mcP_h(\Sigma)} (\langle \bft \cdot \bfsig_{B,\Sigma}(\nablas v) \rangle, 
[\nablas w])_\bfx
\\ \nonumber
&\qquad 
-\sum_{\bfx\in\mcP_h(\Sigma)} ([\nablas v] \textcolor{black}{,} \langle \bft \cdot \bfsig_{B,\Sigma}(\nablas w), \rangle 
)_\bfx
 \\ \nonumber
&\qquad  
+ \sum_{\bfx\in\mcP_h(\Sigma)} \beta_{\Sigma,z} ([\nablas v ], [\nablas w])_{\bfx}
\\
s_{\Sigma,h}(v,w) &= \sum_{F\in \mcF_{h,I}(\Sigma)} \sum_{j =1}^k \gamma_{\Sigma,1} h^{2(j - 2)} ([\partial^j_{\bfn_F} v ] ,[\partial^j_{\bfn_F} w ])_F  
\\ \nonumber
&\qquad + \sum_{T \in \mcT_{h,I}(\Sigma)} \sum_{j=0}^2 \gamma_{\Sigma,2} h^{2(j - 2)+1} 
  (\partial_{\bfn_\Sigma} \partial^j_{\bft} v, \partial_{\bfn_\Sigma} \partial^j_{\bft} w )_T
\end{align}
%
the penalty parameter takes the form
\begin{equation}
\beta_{\Sigma} = \beta_{\Sigma,0} \mcCB =  \beta_{\Sigma,0}  { \Es I_\Sigma}
\end{equation}
with $\beta_{\Sigma,0}$ a parameter that only depends on the polynomial order, and 
$s_h$, with positive parameters $\gamma_{\Sigma,i}$, is a stabilization term which 
is added to ensure coercivity and stability of the stiffness matrix, cf. \cite{BurClaHan15}. 
\begin{rmk} Using the identities 
\begin{equation}
\nablas v = (\partial_{\bft} v ) \bft,
\qquad \bfeps_\Sigma(\nablas v ) = (\partial_{\bft}^2 v) \bft \otimes \bft,
\qquad 
\bfsig_{B,\Sigma}(\nablas v ) = \Es I_\Sigma (\partial_{\bft}^2 v) \bft \otimes \bft
\end{equation}
we note that $a_{\Sigma,h}$ can alternatively be written 
in the form
\begin{align}
a_{\Sigma,h}(v,w) \label{asigma}
&=
\sum_{T\in \mcT_h(\Sigma)}  
(\Es I_\Sigma  \,\partial_{\bft}^2 v , \partial_{\bft}^2 w)_{\Sigma\cap T} 
\\ \nonumber
&\qquad 
-\sum_{\bfx\in\mcP_h(\Sigma)} 
( \langle \Es I_\Sigma \partial_{\bft}^2  v \rangle, [\partial_{\bft} w])_{\bfx}
\\ \nonumber
 &\qquad  -\sum_{\bfx\in\mcP_h(\Sigma)} 
 ([\partial_{\bft}\textcolor{black}{ v}],\langle \Es I_\Sigma \partial_{\bft}^2 \textcolor{black}{w} \rangle )_{\bfx}
 \\ \nonumber
&\qquad  + \sum_{\bfx\in\mcP_h(\Sigma)} \frac{\beta_{\Sigma,0}}{h} (\Es I_\Sigma[\partial_{\bft} v ],[\partial_{\bft} w] )_{\bfx}
\end{align}
which is the form in \cite{EngGar02}.
\end{rmk}
\begin{rmk}
The terms on the discrete set $\mcP_h(\Sigma)$ are associated with the work 
of the end moments on the end rotation which occur due to the lack of $C^1(\Omega)$ 
continuity of the approximation, as in the plate model. See Remark \ref{rem:bending-moments}.
\end{rmk}
\begin{rmk} We note that due to the stabilization this method works for a single 
beam, i.e. without being embedded in a plate. The basic principle is the same as 
for the trace finite element method proposed in \cite{OlReGr09} and the stabilized 
version proposed in \cite{BuHaLa15}. \textcolor{black}{
When the beam is embedded in a plate, which is the case in this work, the need for the stabilization term is 
mitigated, and if the plate is sufficiently stiff we may omit the stabilization term, see Section \ref{sec:coercivity} for further details.
}
\end{rmk}

\subsection{The c/dG Method for the Reinforced Plate Model}

Recall that $\mcS = \{ \Sigma \}$ is a set of beams arbitrarily oriented in $\Omega$. 
Using superposition we obtain the problem: find $u_h \in V_{\Omega,h}$ such that 
\begin{equation}
A_h(u_h, v) = l(v)\qquad \forall v \in V_{\Omega,h}
\end{equation}
where  the forms are defined by
\begin{align}
A_h(v,w) & = a_{\Omega,h}(v,w) + \sum_{\Sigma \in \mcS} a_{\Sigma,h}(v,w)
\\
l(v) &= l_{\Omega} (v) + \sum_{\Sigma \in \mcS} l_{\Sigma}(v)
\end{align}

\subsection{Coercivity for Reinforced Plates}
\label{sec:coercivity}
In this section we study the coercivity of the c/dG method for the reinforced plate. We shall use the
stability provided by the plate to prove stability of the reinforced
plate, \textcolor{black}{without the need of the stabilizing
  terms ($\gamma_{\Sigma,1} = \gamma_{\Sigma,2} = 0$)}. \textcolor{black}{This is only possible as long as the mesh
  size $h$ is larger than or equal to the beam with $b_{\Sigma}$. When
  this condition is not satisfied, stability uniform in $h$ is achieved only when the
  stabilizing terms are included  ($\gamma_{\Sigma,1} ,
  \gamma_{\Sigma,2} > 0$), using similar ideas as in \cite{BuHaLa15, BHLM16}.}

\paragraph{Coercivity of the Plate} We first recall that 
the c/dG method for the plate is coercive. Introducing the energy norm
\begin{equation}
\tn v \tn^2_{\Omega,h} = 
 \sum_{T \in \mcT_{h}} 
\mcCP  \| \nabla^2 v  \|^2_{T} 
+ 
\sum_{F \in \mcF_h}  \mcCP  h \| \langle \nabla^2 v \rangle \|^2_{F}
+
\sum_{F \in \mcF_h} \mcCP  h^{-1} \| [\nabla v] \|^2_{F}
\end{equation}
there is a constant $m_P>0$ such that 
\begin{equation}\label{eq:coercivity-plate}
m_P \tn v \tn^2_{\Omega,h} \leq a_{\Omega,h} (v,v) \qquad \forall v \in V_{\Omega,h} 
\end{equation}
for $\beta_{\Omega}$ large enough.

\paragraph{Coercivity of the Reinforced Plate}
Next turning to the reinforced plate we introduce the energy norm associated with 
the beam
\begin{equation}
\tn v \tn^2_{\Sigma,h} = 
  \sum_{T \in \mcT_{h}(\Sigma)} 
\mcCB \| \textcolor{black}{\partial_{\bft}^2} v  \|^2_{ \textcolor{black}{\Sigma \cap T}} 
+ 
\sum_{\bfx \in \mcP_h(\Sigma)} \mcCB h \| \langle \textcolor{black}{\partial_{\bft}^2}  v \rangle \|^2_{\bfx}
+
\sum_{\bfx \in \mcP_h(\Sigma)} \mcCB h^{-1} \| [ \textcolor{black}{\partial_{\bft}}   v] \|^2_{\bfx}
\end{equation}
Then there is a constant $m$ such that 
\begin{equation}\label{eq:coercivity-reinforced}
m \left( \tn v \tn^2_{\Sigma,h} + \tn v \tn^2_{\Omega,h} \right) 
\lesssim A_h(v,v)\qquad \forall v \in V_h
\end{equation}
for $\beta_\Omega$ and $\beta_\Sigma$ large enough.
\paragraph{Verification of (\ref{eq:coercivity-reinforced})} Using the following two 
inequalities, which we verify below, 
\begin{equation}\label{eq:step2-main-est}
C_1\left(  \sum_{\bfx \in \mcP_h(\Sigma)} \mcCB h \| \langle  \textcolor{black}{\partial_{\bft}^2} v \rangle \|^2_{\bfx}    \right) 
\leq  \sum_{T \in \mcT_{h}} 
\mcCP  \| \nabla^2 v  \|^2_{T} 
\end{equation}
for some constant $C_1>0$, and 
\begin{align}\label{eq:step1-main-est}
\sum_{T \in \mcT_{h}(\Sigma)} \mcCB  \| \partial^2_{\bft} v \|^2_{\Sigma \cap T} 
 +  \sum_{\bfx \in \mcP_h(\Sigma)} \mcCB h^{-1} \| [ \partial_{\bft} v]\|^2_{\bfx}
 \leq   \sum_{T \in \mcT_{h}} 
\frac{m_P}{3} \mcCP  \| \nabla^2 v  \|^2_{T} + a_{\Sigma,h}(v,v)
\end{align}
for $\beta_\Sigma$ large enough,  we have
\begin{align}
A_h(v,v) &= a_{\Omega,h}(v,v) + a_{\Sigma,h}(v,v) 
\\
&\geq m_P \tn v \tn^2_{\Omega,h} + a_{\Sigma,h}(v,v) 
\\
&= \frac{m_P}{3} \tn v \tn^2_{\Omega,h} + 
 \frac{m_P}{3} \tn v \tn^2_{\Omega,h}  +
\left( \frac{m_P}{3} \tn v \tn^2_{\Omega,h} 
+ a_{\Sigma,h}(v,v)\right) 
\\
&\geq 
 \frac{m_P}{3} \tn v \tn^2_{\Omega,h} 
 + \frac{C_1 m_P}{3}\textcolor{black}{\left(  \sum_{\bfx \in \mcP_h(\Sigma)} \mcCB h \| \langle  \partial_{\bft}^2 v \rangle \|^2_{\bfx}    \right) }
\\
\nonumber 
&\qquad +  
\left(  \sum_{T \in \mcT_{h}(\Sigma)} \mcCB  \| \partial^2_{\bft} v \|^2_{\Sigma \cap T} 
 +    \sum_{\bfx \in \mcP_h(\Sigma)} \mcCB h^{-1} \| [ \partial_{\bft} v]\|^2_{\bfx}\right)
 \\
&\geq m  \left( \tn v \tn^2_{\Sigma,h} + \tn v \tn^2_{\Omega,h} \right) 
\end{align}
where $m = \min(m_P/3,C_1 m_P/3,1)$.

\paragraph{Verification of (\ref{eq:step2-main-est})} 
We note that, for $\bfx \in \Sigma \cap T$, $T\in \mcT_{h}$, we have the inverse inequality  
\begin{align}\label{eq:inverse-coer}
\| \partial_{\bft}^2 v\|_{\bfx} 
&\leq C_{\text{inv}} h^{-1} \|\partial_{\bft}^2 v \|_{T} 
\leq C_{\text{inv}}  h^{-1}  \|\nabla^2 v \|_T 
\end{align}
Using (\ref{eq:inverse-coer}) we obtain, \textcolor{black}{with $\mcT_{h}(\bfx) = \{ T \in \mcT_h : \bfx \in \overline{T} \}$},
\begin{align}
\sum_{\bfx \in \mcP_h(\Sigma)} \mcCB h \| \langle\textcolor{black}{\partial_{\bft}^2}  v\rangle \|^2_{\bfx} 
&\leq
\sum_{\bfx \in \mcP_h(\Sigma)}  C^2_{\text{inv}} \frac{\mcCB}{ \mcCP h} 
\mcCP \| \nabla^2 v \|^2_{\mcT_h(\bfx)} 
\\
&\leq C^2_{\text{inv}} \frac{\mcCB}{ \mcCP h} \tn v \tn^2_{\textcolor{black}{\Omega,h}}
\end{align}
and thus we have the estimate 
\begin{align}
\underbrace{\frac{1}{C^2_{\text{inv}}} \frac{\mcCP h}{\mcCB}}_{C_1}
\left( \sum_{\bfx \in \mcP_h(\Sigma)} \mcCB h \| \langle\textcolor{black}{\partial_{\bft}^2}  v \rangle \|^2_{\bfx} \right)
\leq  \tn v \tn^2_{\textcolor{black}{\Omega,h}}  
\end{align}
We note, using the definitions (\ref{eq:mcCP}) and (\ref{eq:mcCB}) of $\mcCP$ and 
$\mcCB$, that 
\begin{equation}\label{eq:quotient-constants}
 \frac{\mcCP h}{\mcCB } 
 =  \frac{1}{1+\nub} \frac{\Eb t_\Omega^3}{ \Es t_\Sigma^3}  
 \frac{h}{b_\Sigma} 
 \geq 
 \frac{1}{1+\nub} \frac{\Eb t_\Omega^3}{ \Es t_\Sigma^3}  
 C_{\text{mesh}} 
\end{equation}
where we used \textcolor{black}{the condition that the beam width is
  smaller than the mesh size} (\ref{eq:Cmesh}) and thus the right hand side is a positive constant 
independent of the mesh size and so is $C_1$.

\paragraph{Verification of (\ref{eq:step1-main-est})}   First we have the \textcolor{black}{equality}
\begin{align}
a_{\Sigma,h}(v,v) 
&\textcolor{black}{= }
\sum_{T \in \mcT_{h}(\Sigma)} \mcCB  \| \partial^2_{\bft} v \|^2_{\Sigma \cap T} 
\\ \nonumber
&\qquad -
\textcolor{black}{2} \underbrace{\sum_{\bfx \in \mcP_h(\Sigma)} \mcCB ( \langle  \partial^2_{\bft} v \rangle, [\partial_{\bft} v ])_\bfx }_{\bigstar}
\\ \nonumber 
&\qquad
+
\sum_{\bfx \in \mcP_h(\Sigma)} \beta_{\Sigma,0} \mcCB h^{-1} \| [\partial_{\bft} v] \|^2_{\bfx} 
\end{align}
To estimate $\bigstar$  we employ the inverse inequality (\ref{eq:inverse-coer}) as follows
\begin{align}
\bigstar &= \textcolor{black}{2} \sum_{\bfx \in \mcP_h(\Sigma)} \mcCB ( \langle  \partial^2_{\bft} v \rangle, [\partial_{\bft} v ])_\bfx
\\
&\leq
\textcolor{black}{2} \sum_{\bfx \in \mcP_h(\Sigma)}
\mcCB  C_{\text{inv}}  h^{-1}  \| \nabla^2 v \|_{\mcT_h(\bfx)} \| [\partial_{\bft} v ]\|_{\bfx}
\\ 
&\leq
\sum_{T \in \mcT_{h}(\Sigma)} 
\delta \mcCB  C_{\text{inv}}^2  h^{-1} \| \nabla^2 v \|^2_{T} 
\\ \nonumber
&\qquad + 
\sum_{\bfx \in \mcP_h(\Sigma)} \delta^{-1} \mcCB  h^{-1} \| [\partial_{\bft} v] \|^2_{\bfx}
\end{align}
where we used the inequality $ab \leq (\delta a^2 + \delta^{-1}
b^2)/2$ for $\delta>0$.
We then obtain (\ref{eq:step1-main-est}) 
as follows
\begin{align}
\sum_{T \in \mcT_{h}} 
\frac{m_P}{3}\mcCP  \| \nabla^2 v  \|^2_{T} + a_{\Sigma,h}(v,v)
&\geq 
\sum_{T \in \mcT_{h}(\Sigma)} \mcCB  \| \partial^2_{\bft} v \|^2_{\Sigma \cap T} 
\\ \nonumber
&\qquad +  \sum_{T \in \mcT_{h}} 
\underbrace{(\frac{m_P}{3} \mcCP - \delta \mcCB  C_{\text{inv}}^2  h^{-1} )}_{\geq 0}  
\| \nabla^2 v  \|^2_{T} 
\\ \nonumber
&\qquad +  \sum_{\bfx \in \mcP_h(\Sigma)} \underbrace{( \beta_{\Sigma,0} - \delta^{-1}  )}_{\geq 1} 
\mcCB h^{-1} \| [\partial_{\bft} v] \|^2_{\bfx} 
\\
&\geq \sum_{\bfx \in \mcP_h(\Sigma)} 
\mcCB h^{-1} \| [\partial_{\bft} v] \|^2_{\bfx} 
+
\sum_{\bfx \in \mcP_h(\Sigma)} \mcCB h^{-1} \| [\partial_{\bft} v] \|^2_{\bfx} 
\end{align}
Here we choose: $\delta$ small enough to guarantee that
\begin{equation}
0\leq 
\frac{m_P}{3}\mcCP - \delta \mcCB C^2_{\text{inv}} h^{-1} 
= 
\frac{m_P}{3}\mcCP  \left( 1  - \delta \frac{3}{m_P} \frac{\mcCB C^2_{\text{inv}} }{ \mcCP h }  \right)
=
\frac{m_P}{3}\mcCP  \left( 1  - \delta \frac{3}{m_P} \frac{1}{ C_1 }  \right)
\end{equation} 
where as above, see (\ref{eq:quotient-constants}), $C_1>0$ independent of the mesh 
parameter $h$, and $\beta_\Sigma$ such that 
\begin{align}
\beta_\Sigma - \frac{1}{\delta} \geq 1
\end{align}

\section{Numerical Examples}

In this Section, we give some elementary examples of what can be achieved with the presented technique.
In all numerical examples we use polynomial order $2$, $f_\Sigma=0$, and 
\[
f_\Omega=8\textcolor{black}{\mcCP}(3(x^2(1-x)^2+y^2(1-y)^2)+(1-6x(1-x))(1-6y(1-y)))
\]
corresponding to the solution 
$u = x^2(1 - x)^2y^2(1 - y)^2$ for a clamped plate unsupported by beams.

\textcolor{black}{In order to handle more general boundary conditions we in particular need to be able to impose end displacements on the beam in the case of a free plate (we note that strongly imposed boundary conditions on the plate are also enforced on the beam). Zero displacement of
the beam endpoints $\bfx_E$ are imposed by adding penalty terms
\begin{equation}
\frac{\tilde{\beta}_{\Sigma,0}}{h^3} (\Es I_\Sigma  v ,  w )_{\bfx_E}
\end{equation}
to the form $a_{\Sigma,h}(v,w)$ in (\ref{asigma}), where $\tilde{\beta}_{\Sigma,0}$ is a penalty parameter. �These terms suffice for optimal order convergence (of the beam approximation) in the case of 
second degree polynomial approximations since the shear forces required
for energy consistency are third derivatives of displacements, and thus equal zero.
}
\subsection{Simply supported plate using beams with different supports}

\begin{figure}
\begin{center}
\includegraphics[scale=0.5]{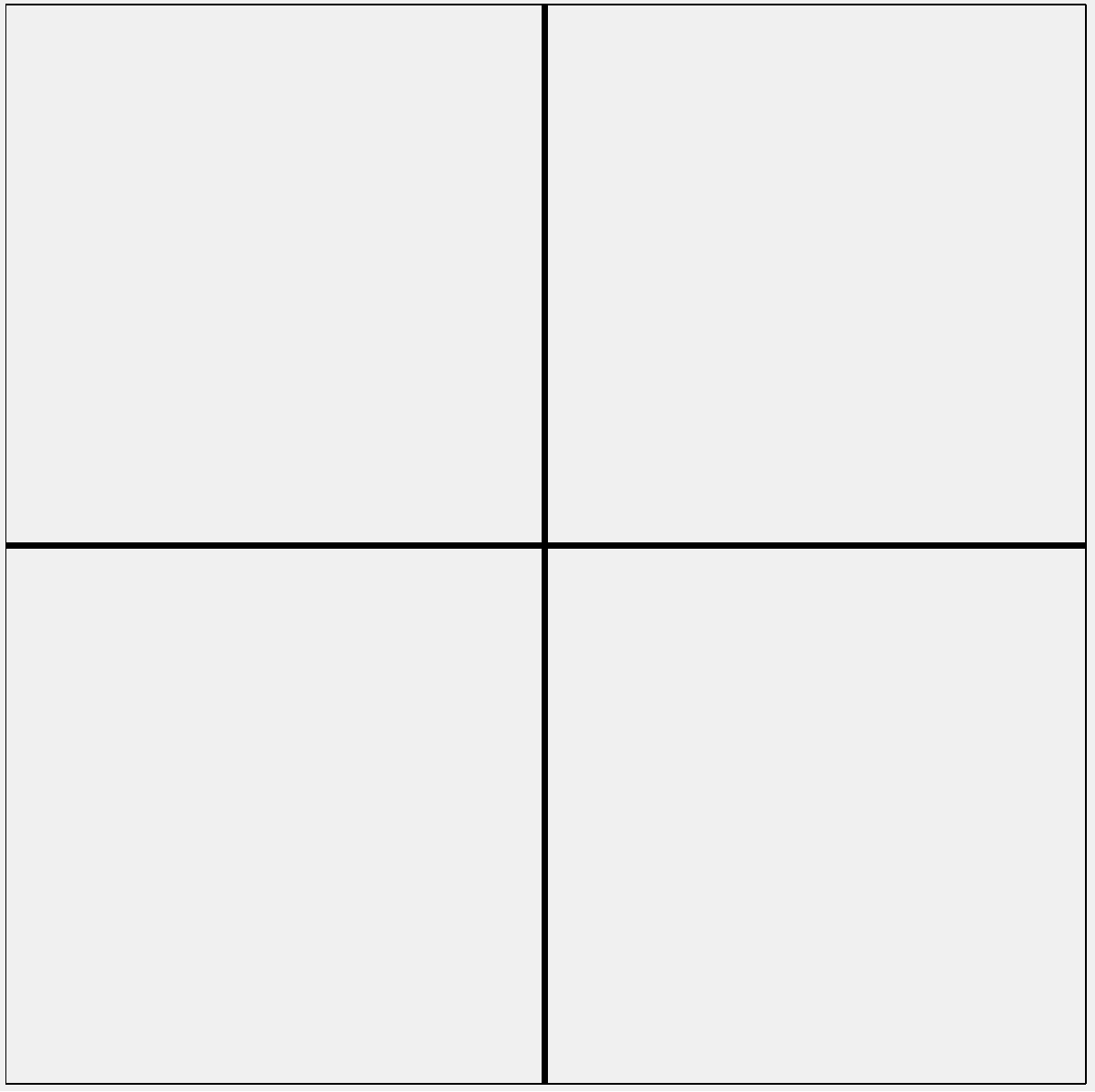}
\end{center}
\caption{Beam reinforced plate.\label{fig:cross}}
\end{figure}

\begin{figure}
\begin{center}
\includegraphics[scale=0.20]{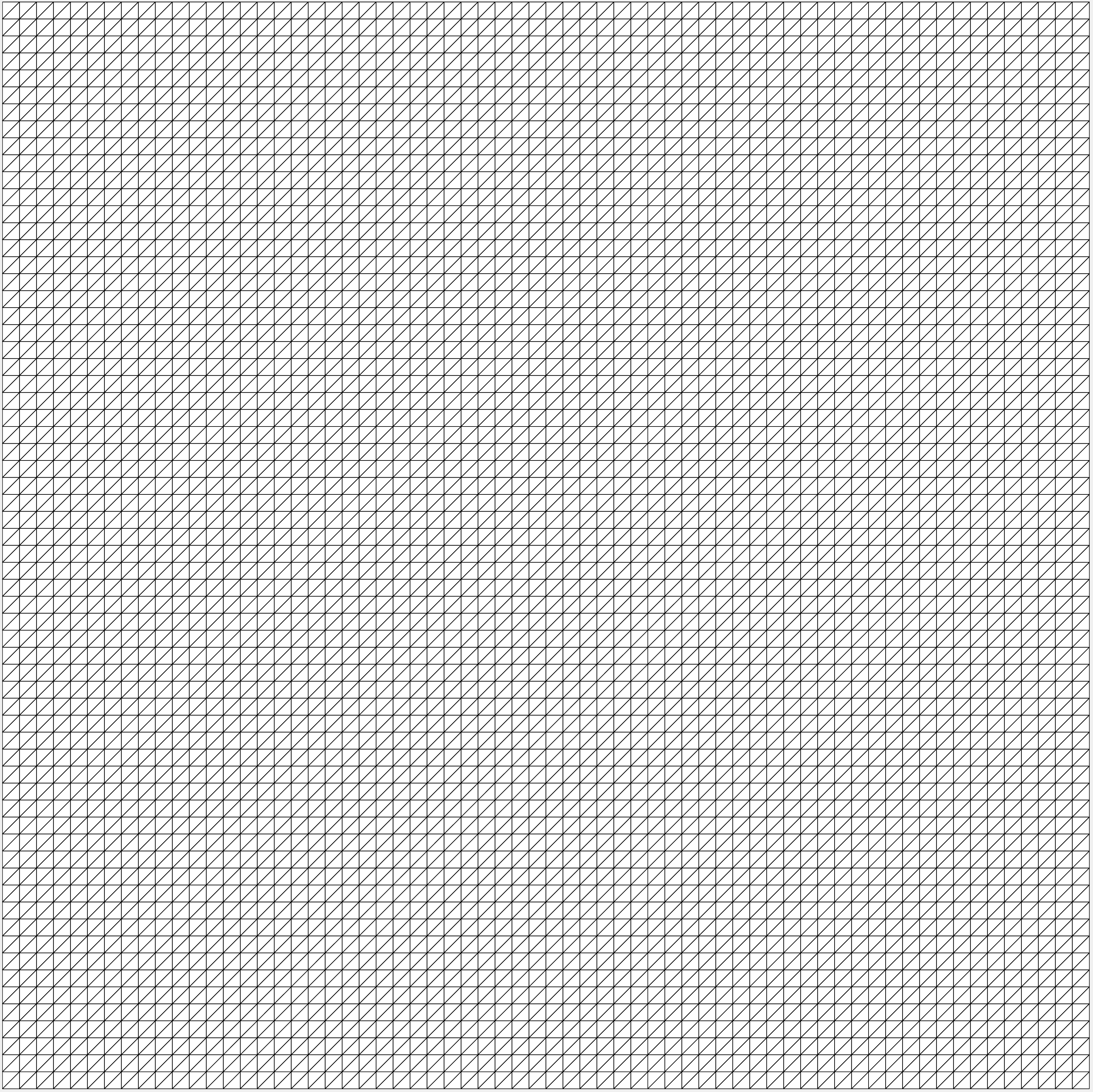}
\end{center}
\caption{Computational mesh.\label{fig:structured}}
\end{figure}

\begin{figure}
\begin{center}
\includegraphics[scale=0.25]{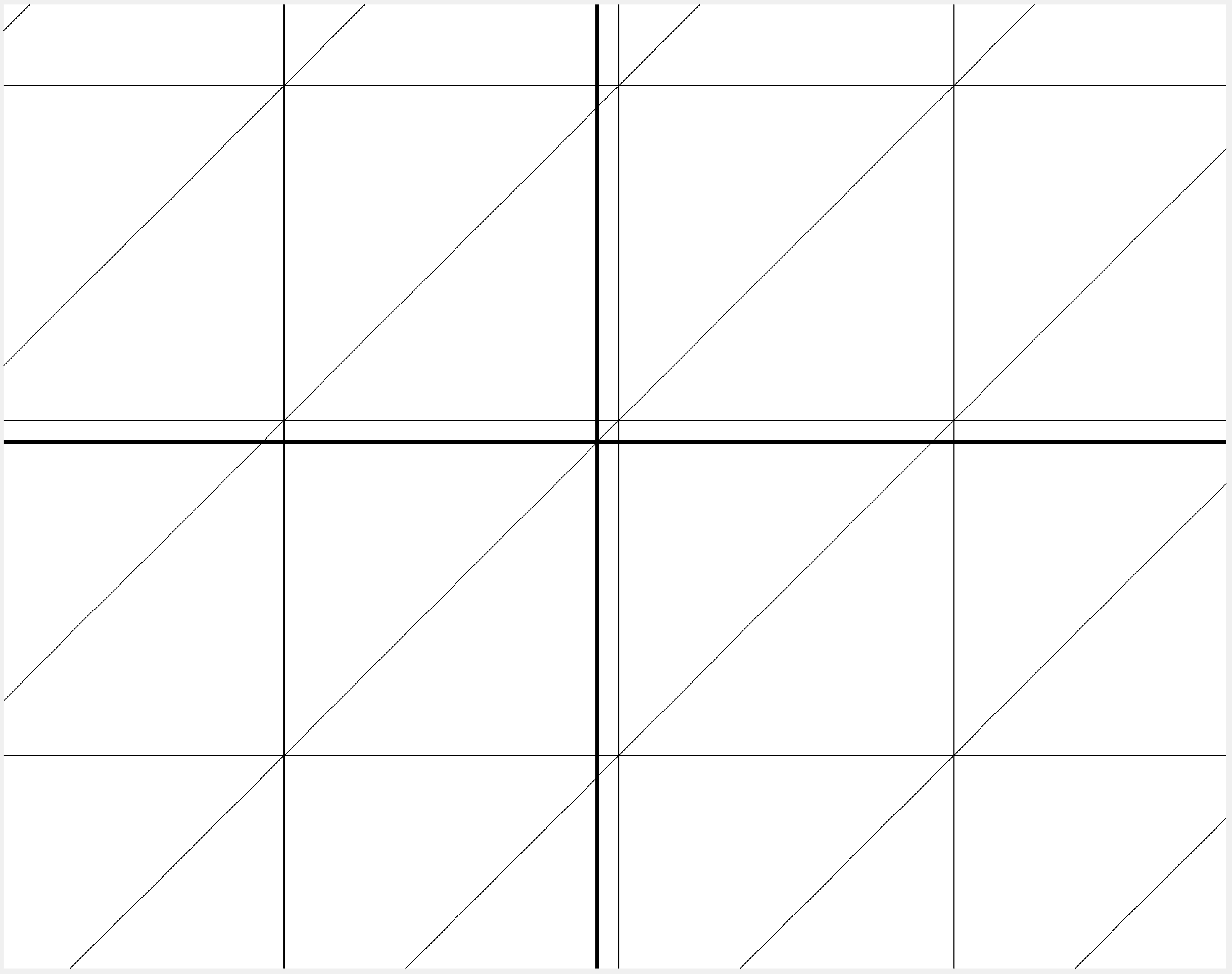}
\end{center}
\caption{Beam/mesh intersection at the center.\label{fig:intersect}}
\end{figure}

We consider a simply supported plate on the domain $\Omega=(0,1)\times(0,1)$ with 
Young's modulus $E_\Omega=100$, Poisson's ratio $\nu_\Omega=1/2$, and thickness $t_\Omega=0.1$.
The plate is supported by two beams oriented as in Fig. \ref{fig:cross},
one at $x=0.499$ and one at at $y=0.499$ (to avoid intersection with the mesh lines).
The computational mesh is shown in Fig. \ref{fig:structured} and in Fig.\ref{fig:intersect} whe
show a close-up of the intersection between the beams and the mesh.

\begin{figure}
\begin{center}
\includegraphics[scale=0.065]{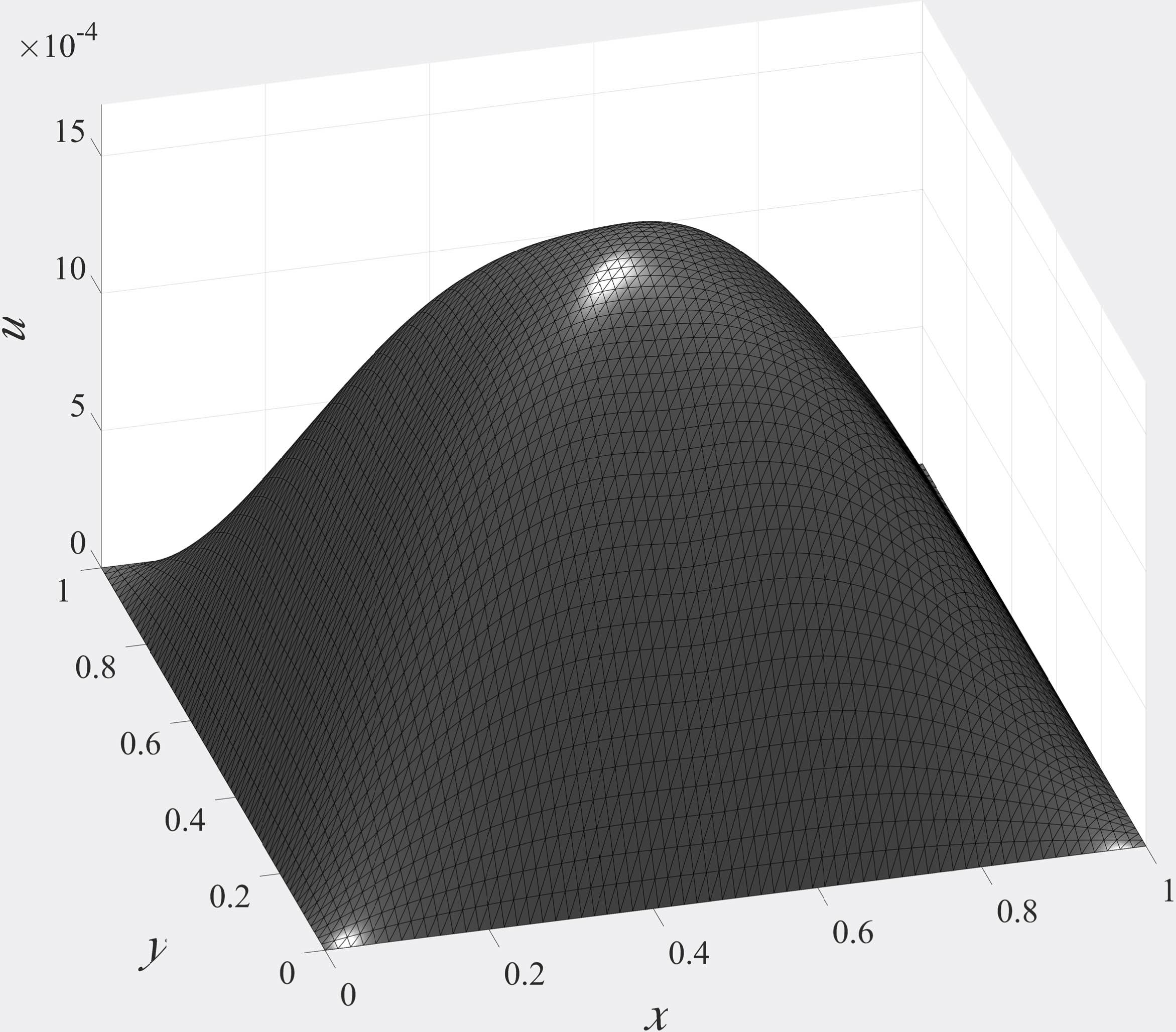}
\includegraphics[scale=0.065]{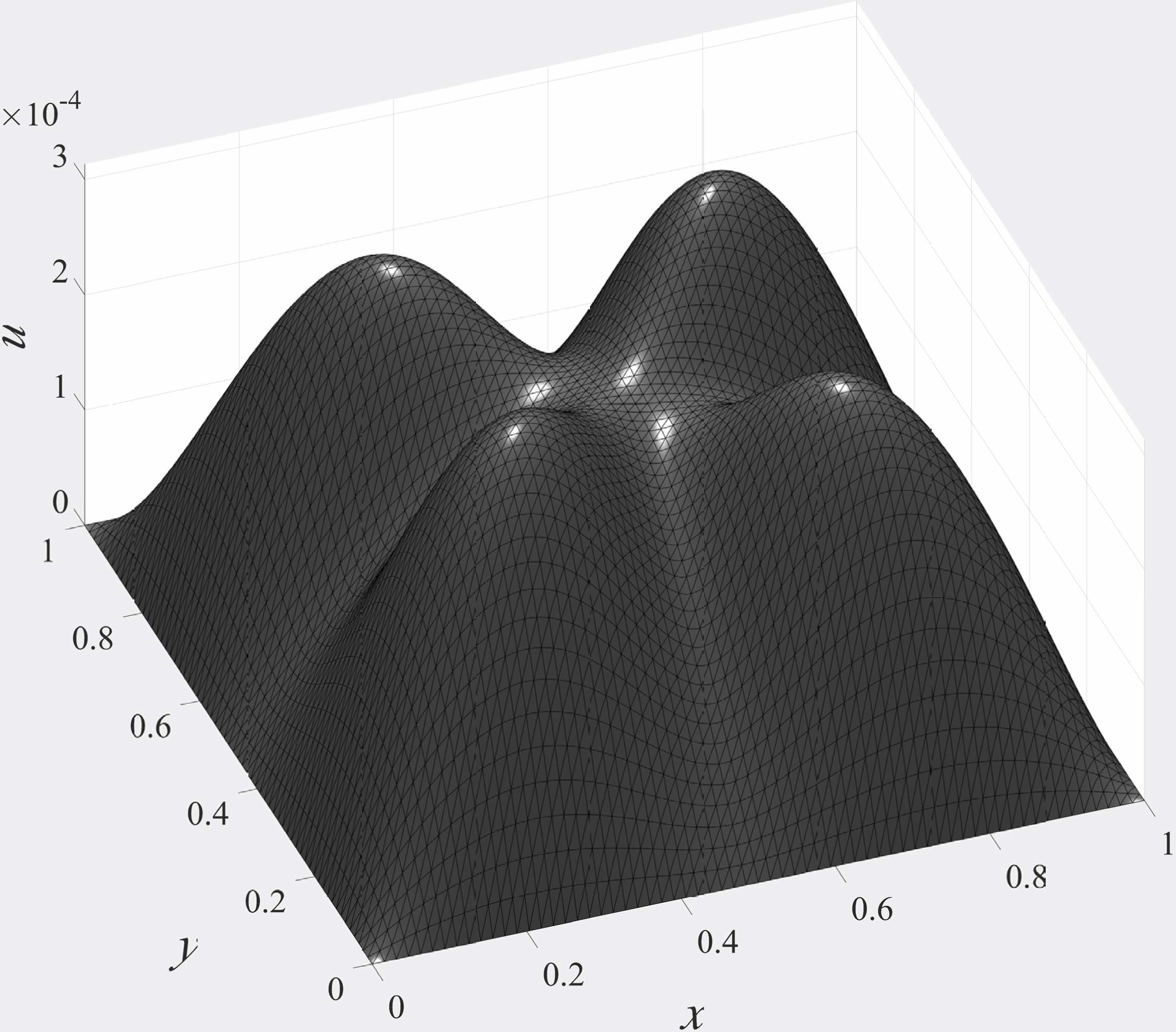}
\end{center}
\caption{Displacements using simply supported support for the beams, $E_\Sigma=100 E_\Omega$ (left) and $E_\Sigma=1000 E_\Omega$ (right).\label{fig:100}}
\end{figure}

\begin{figure}
\begin{center}
\includegraphics[scale=0.35]{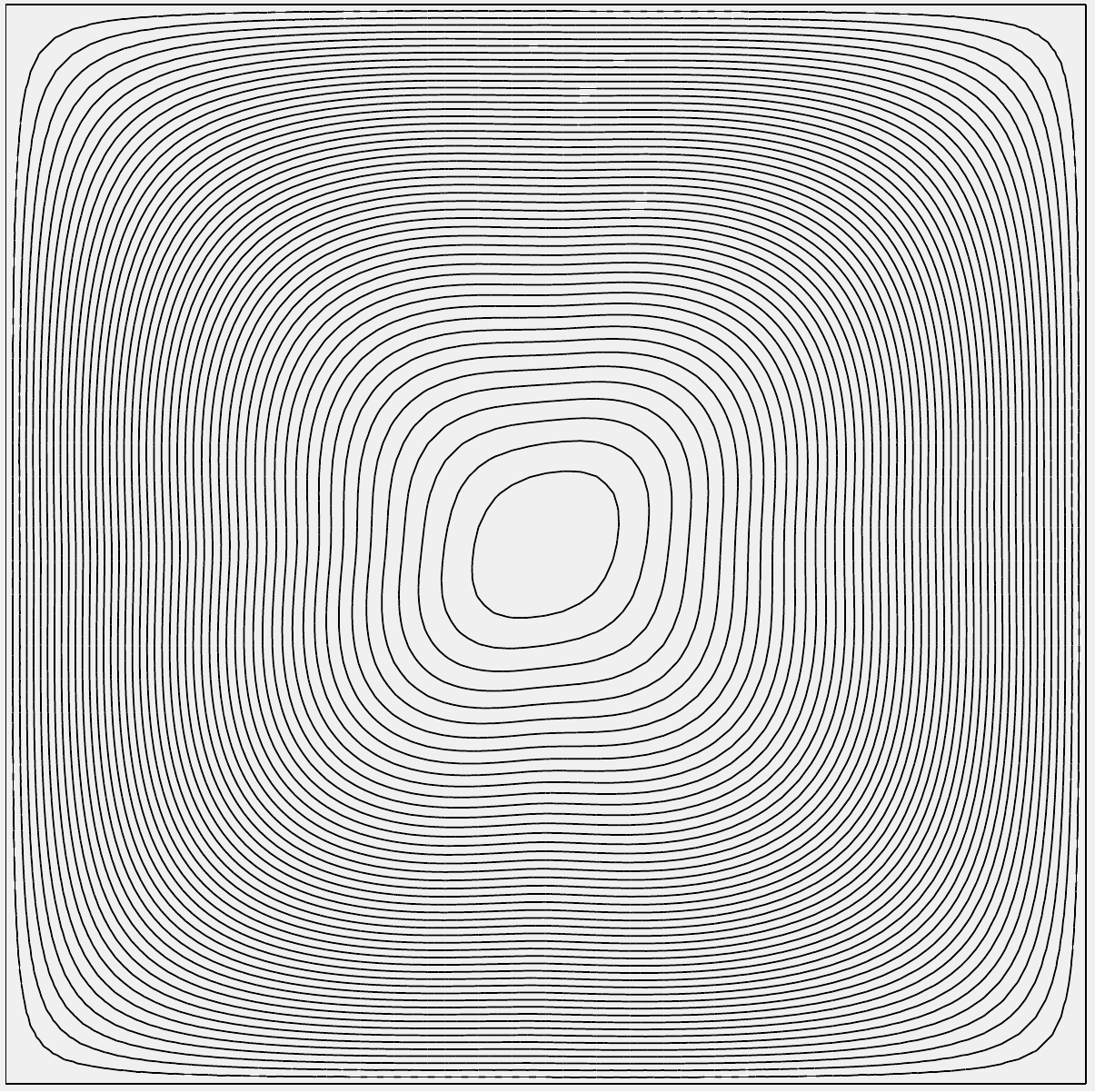}
\includegraphics[scale=0.35]{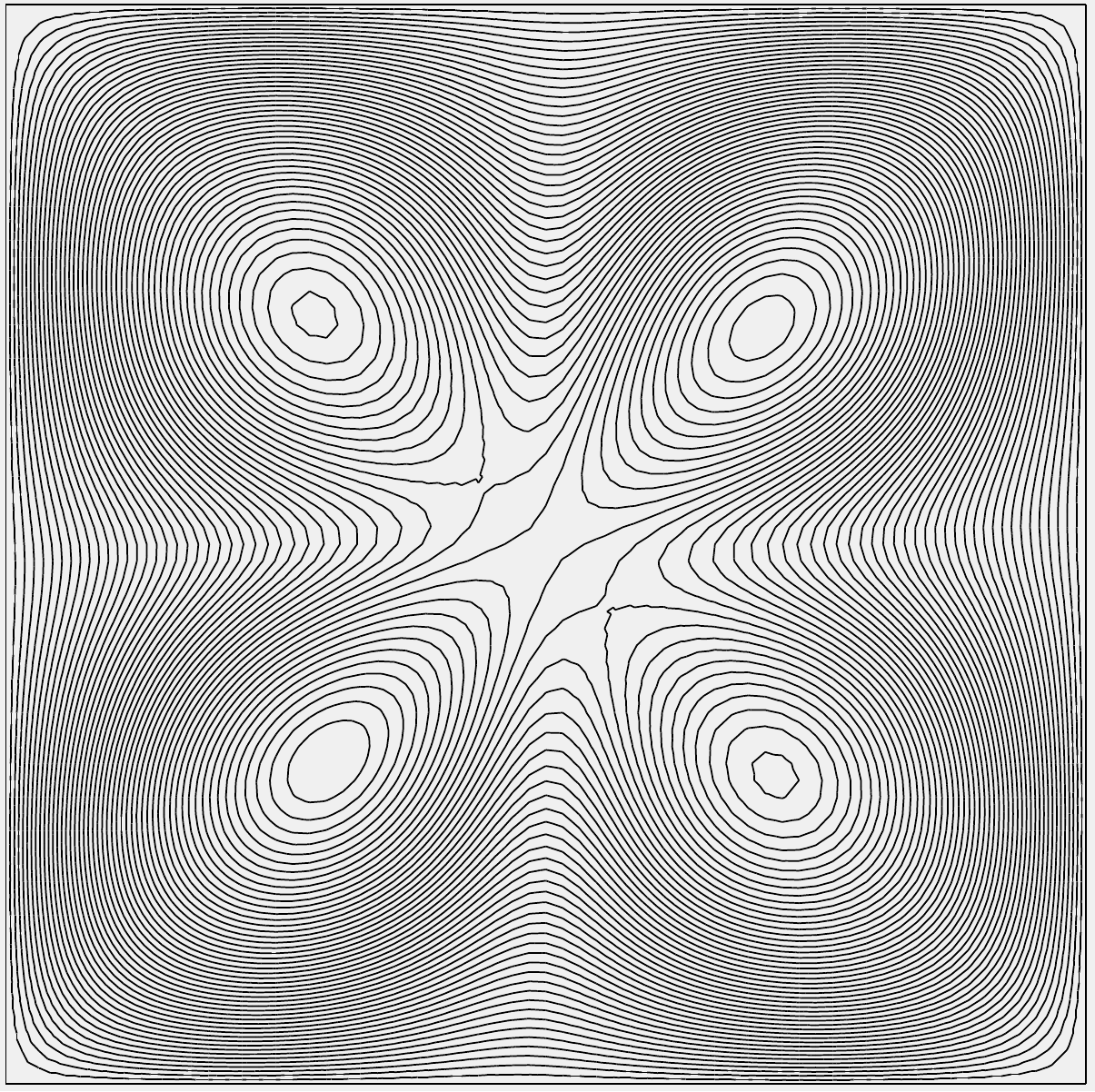}
\end{center}
\caption{Isolines using simply supported beams, $E_\Sigma=100 E_\Omega$ (left) and $E_\Sigma=1000 E_\Omega$ (right).\label{fig:100iso}}
\end{figure}

\begin{figure}
\begin{center}
\includegraphics[scale=0.065]{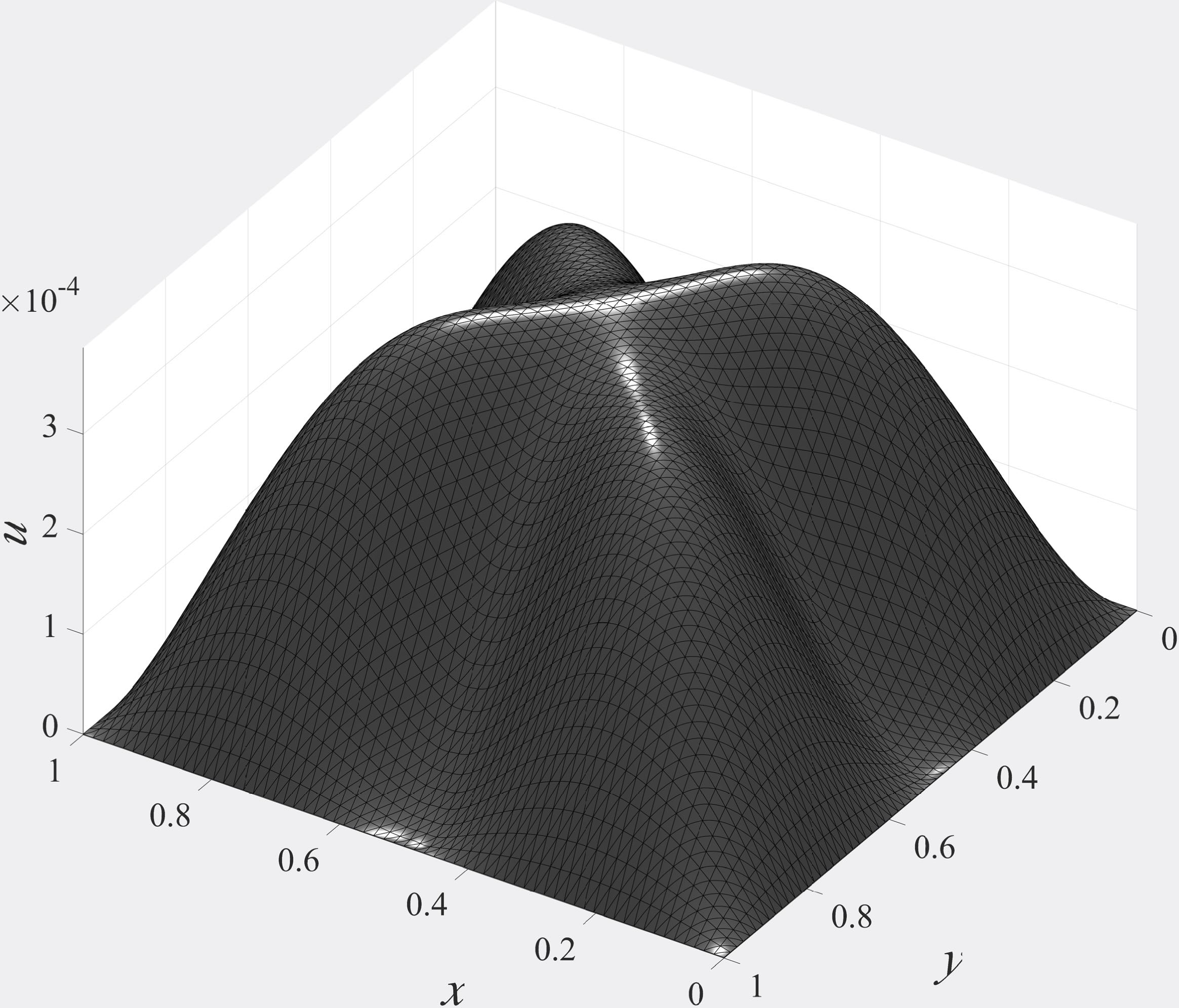}
\includegraphics[scale=0.065]{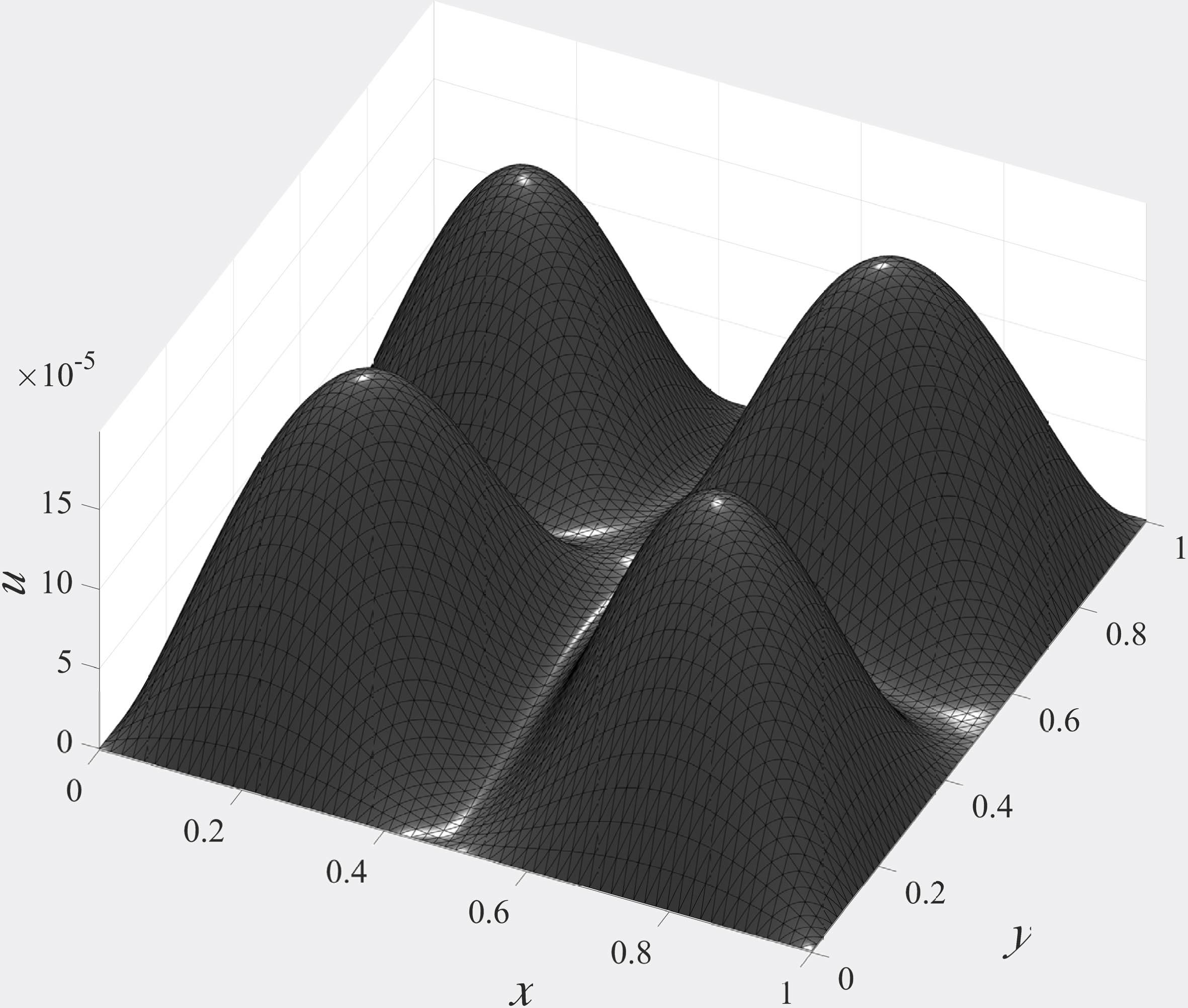}
\end{center}
\caption{Displacements using fixed support for the beams, $E_\Sigma=100 E_\Omega$ (left) and $E_\Sigma=1000 E_\Omega$ (right).\label{fig:1000}}
\end{figure}

\begin{figure}
\begin{center}
\includegraphics[scale=0.35]{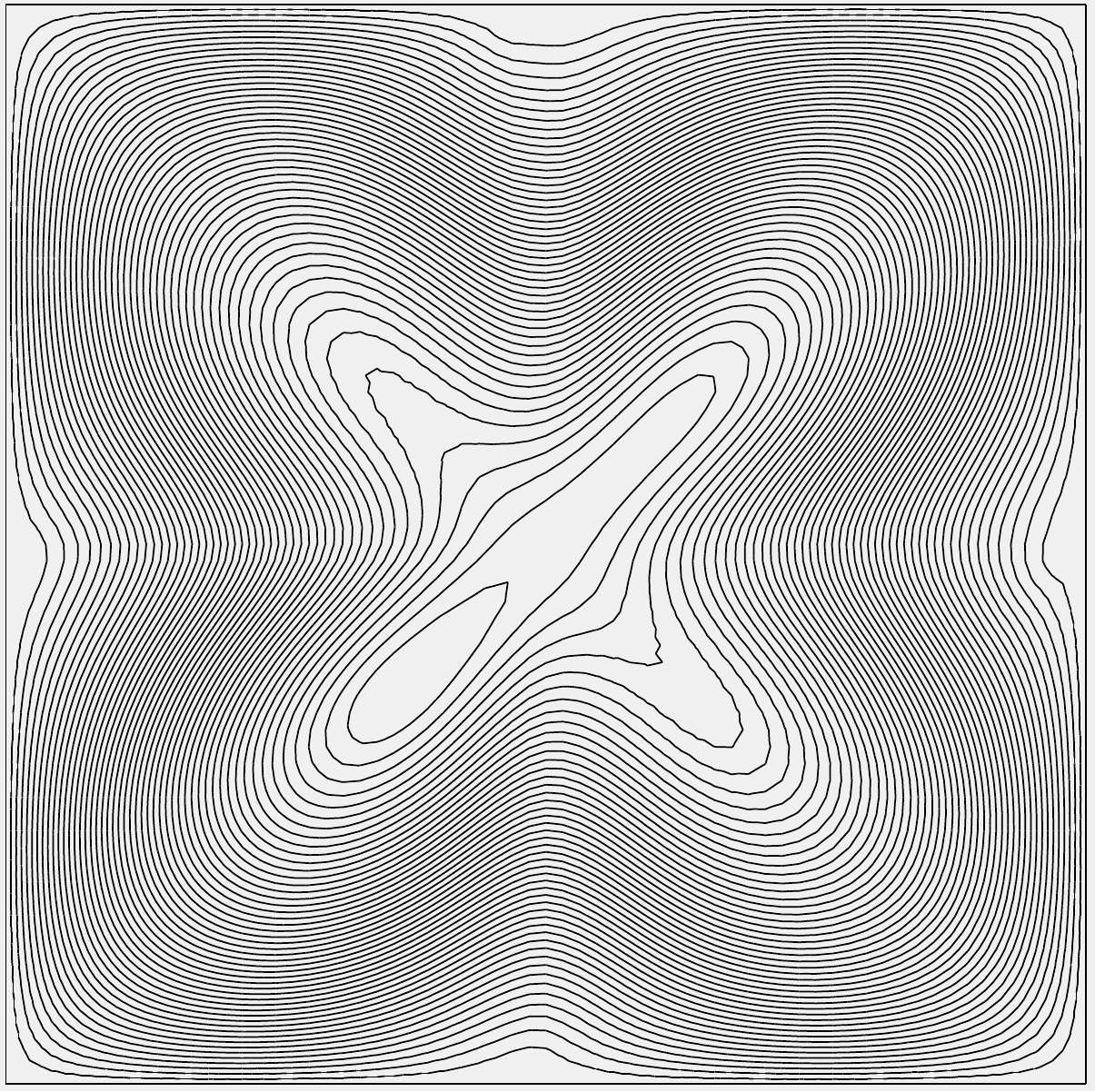}
\includegraphics[scale=0.35]{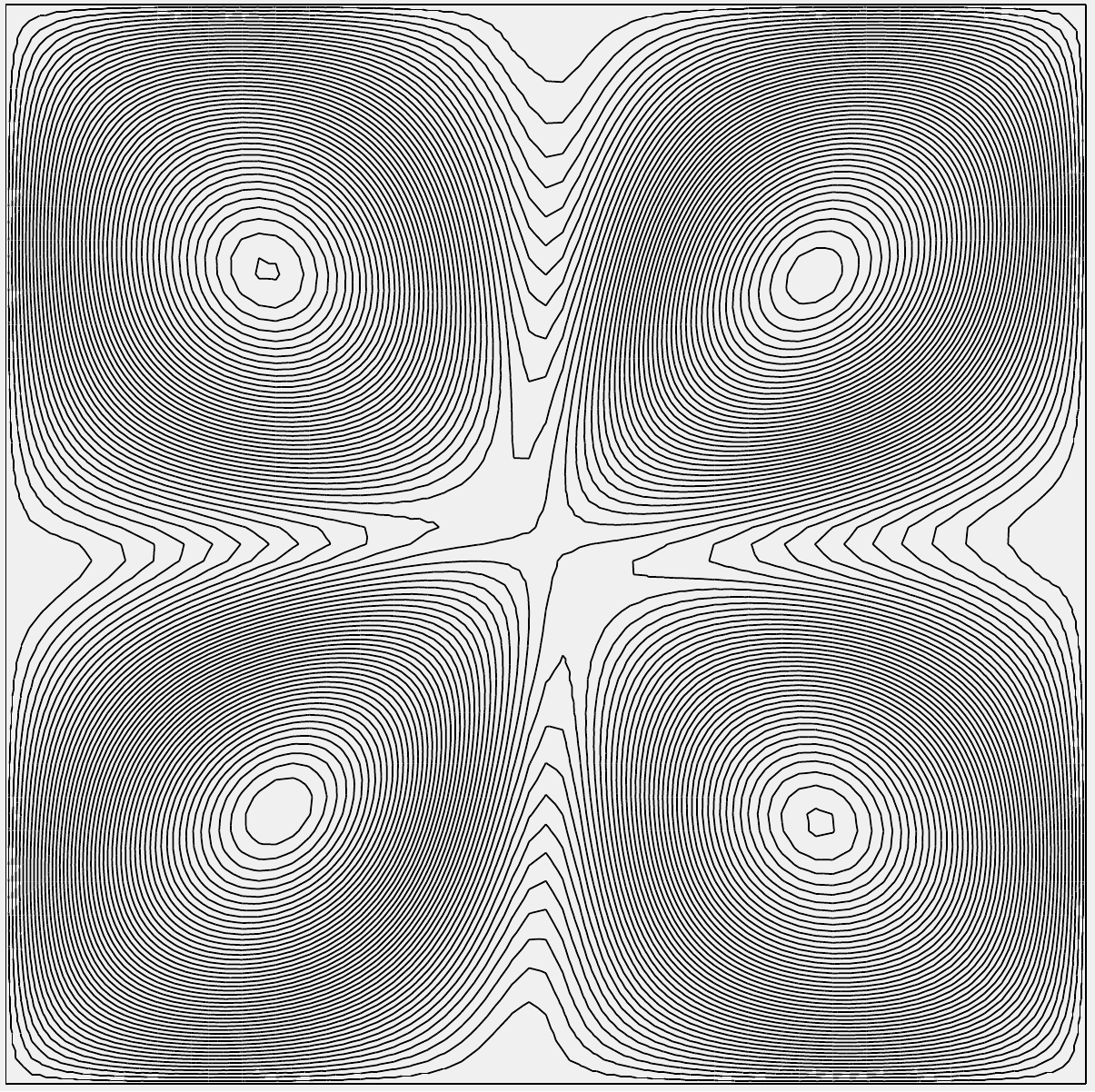}
\end{center}
\caption{Isolines using fixed support for the beams, $E_\Sigma=100 E_\Omega$ (left) and $E_\Sigma=1000 E_\Omega$ (right).\label{fig:1000iso}}
\end{figure}

For this problem we test two different supports for the beams: simply supported and fixed, and two
different stiffnesses for the beams: $E_\Sigma = 100 E_\Omega$ and $E_\Sigma = 1000 E_\Omega$. The thickness and width of the beam are equal and
the same as the thickness of the plate. In Fig. \ref{fig:100} we show the results using $E_\Sigma = 100 E_\Omega$, with simply supported and fixed supports;
in Fig. \ref{fig:100iso} we give the corresponding isolines, and
in Fig. \ref{fig:1000} we show the results using $E_\Sigma = 1000 E_\Omega$, with simply supported and fixed supports;
in Fig. \ref{fig:1000iso} we give the corresponding isolines.

\subsection{Plate only supported by beams}

\begin{figure}
\begin{center}
\includegraphics[scale=0.3]{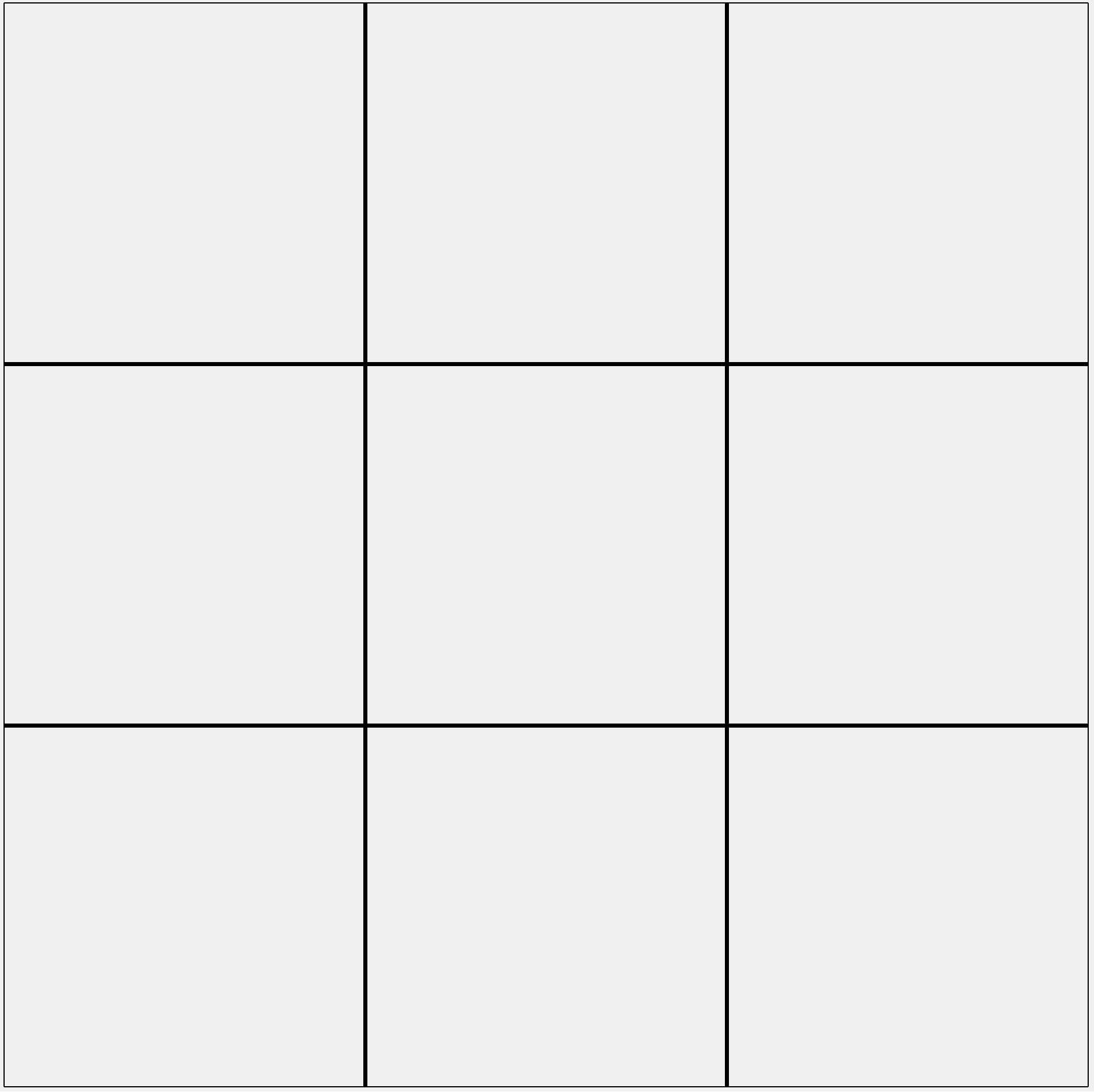}
\end{center}
\caption{Beam reinforced plate.\label{fig:fourbeams}}
\end{figure}

\begin{figure}
\begin{center}
\includegraphics[scale=0.25]{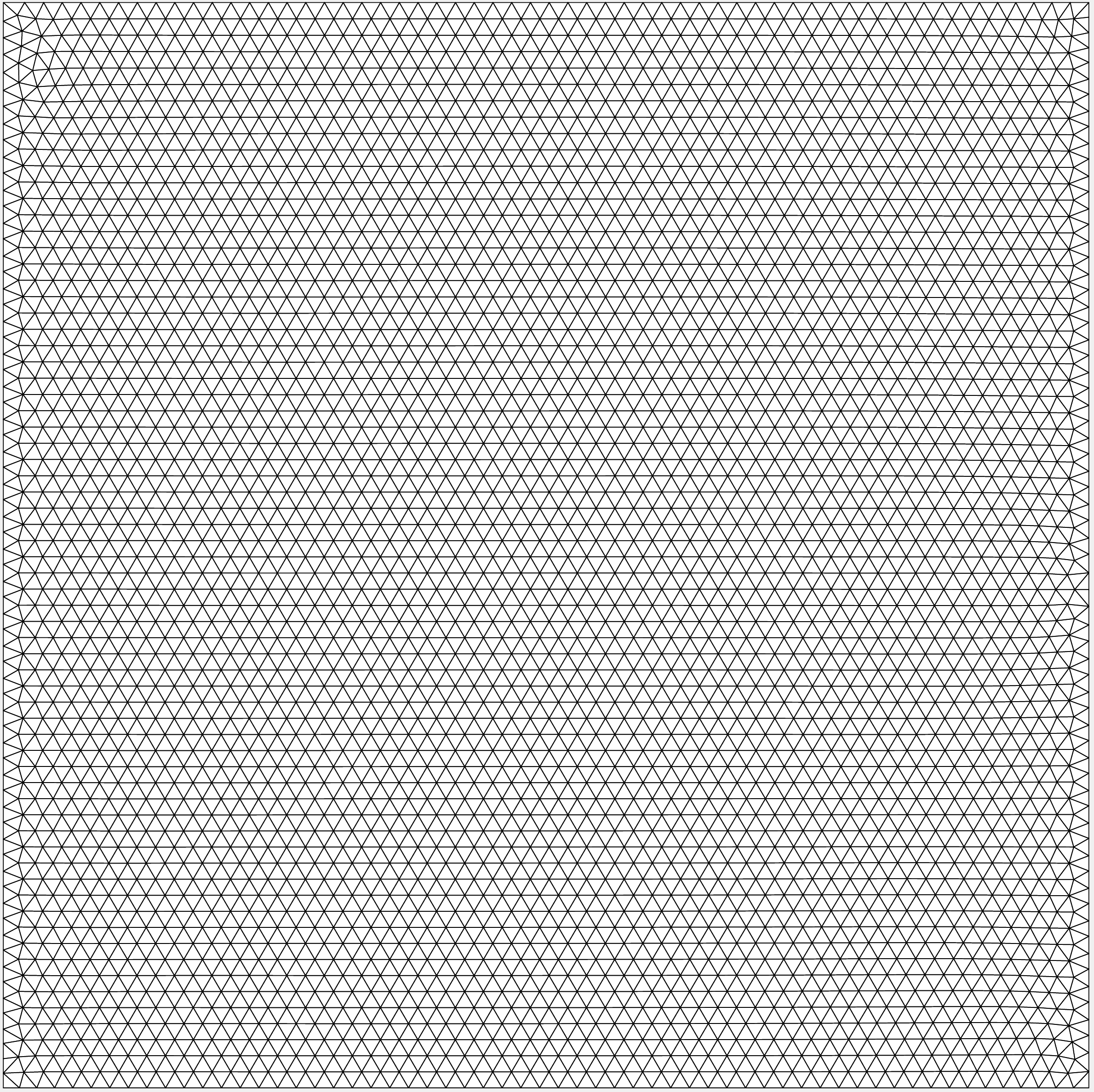}
\end{center}
\caption{Computational mesh.\label{fig:unstruct}}
\end{figure}

Next, we consider a plate with free boundaries, supported only by beams. All data for the plate are the same as in the previous example.
The plate is supported by four beams positioned at $1/3$ and $2/3$ from each boundary as indicated in Fig. \ref{fig:fourbeams}.
The beams have the same dimension as previously, with Young's modulus $E_\Sigma = 100 E_\Omega$. The computational mesh is 
unstructured and shown in Fig. \ref{fig:unstruct}.

\begin{figure}
\begin{center}
\includegraphics[scale=0.065]{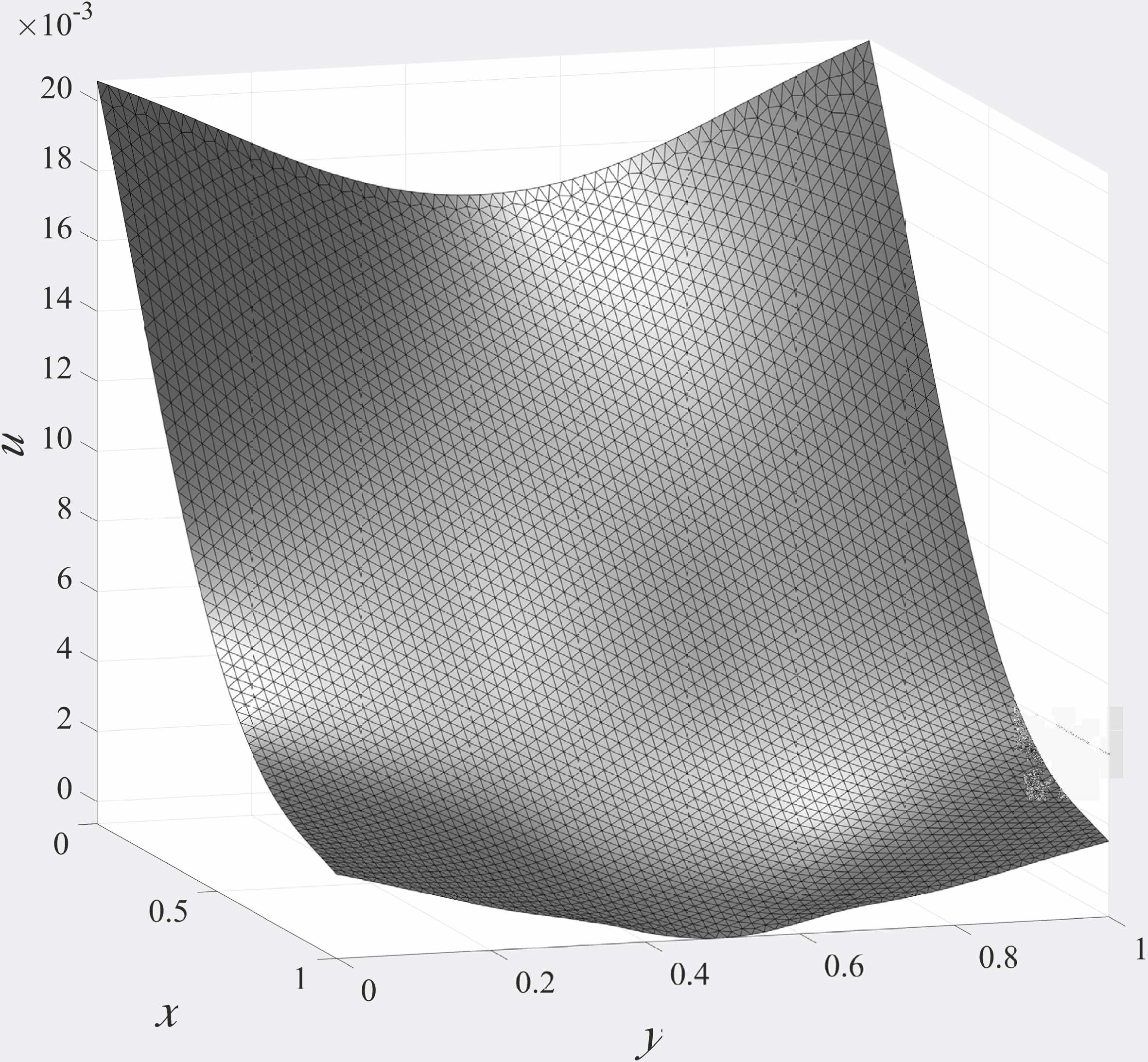}
\includegraphics[scale=0.22]{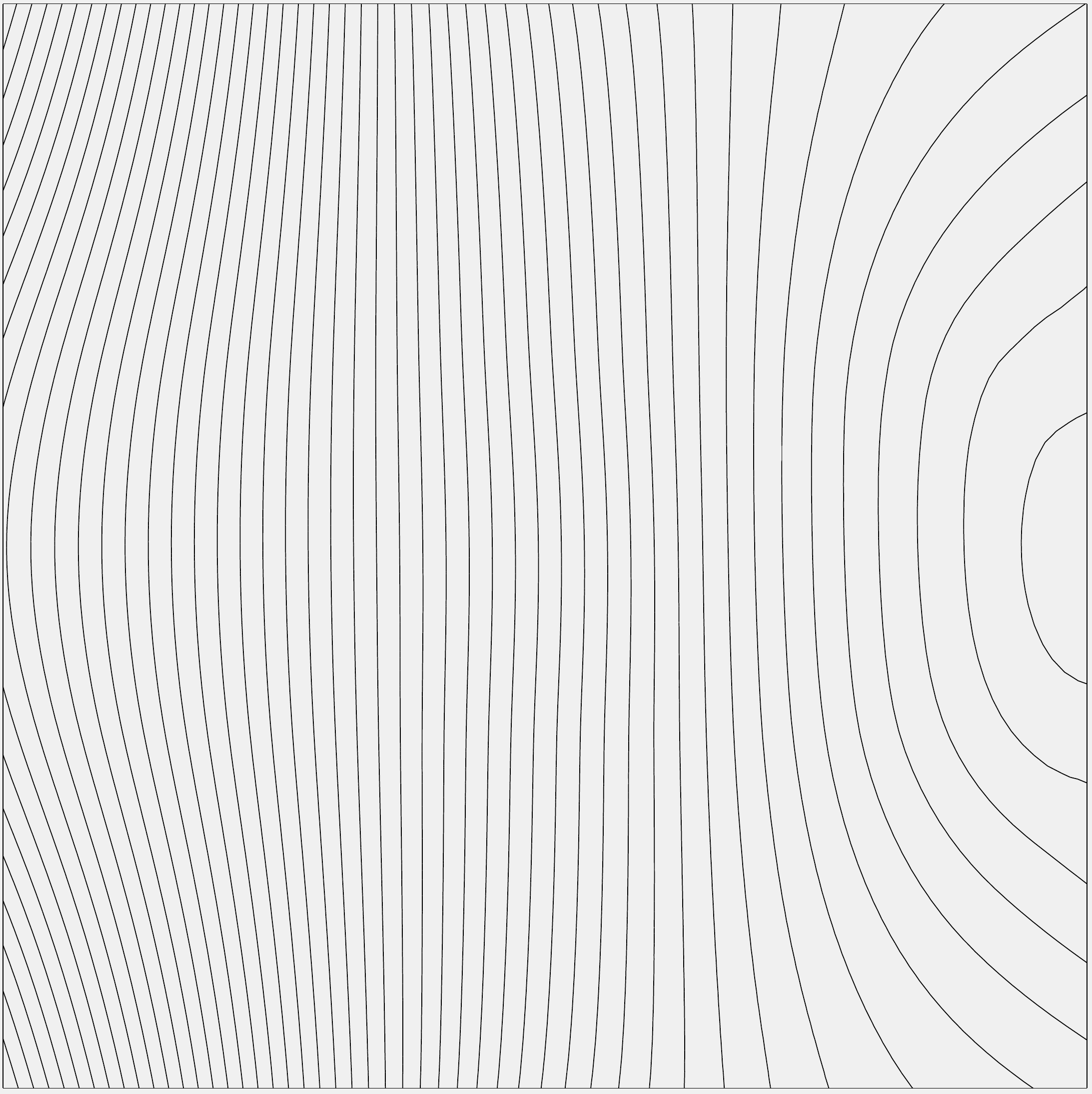}
\end{center}
\caption{Deformations when the beams are clamped at $x=1$.\label{fig:onesidele}}
\end{figure}

\begin{figure}
\begin{center}
\includegraphics[scale=0.15]{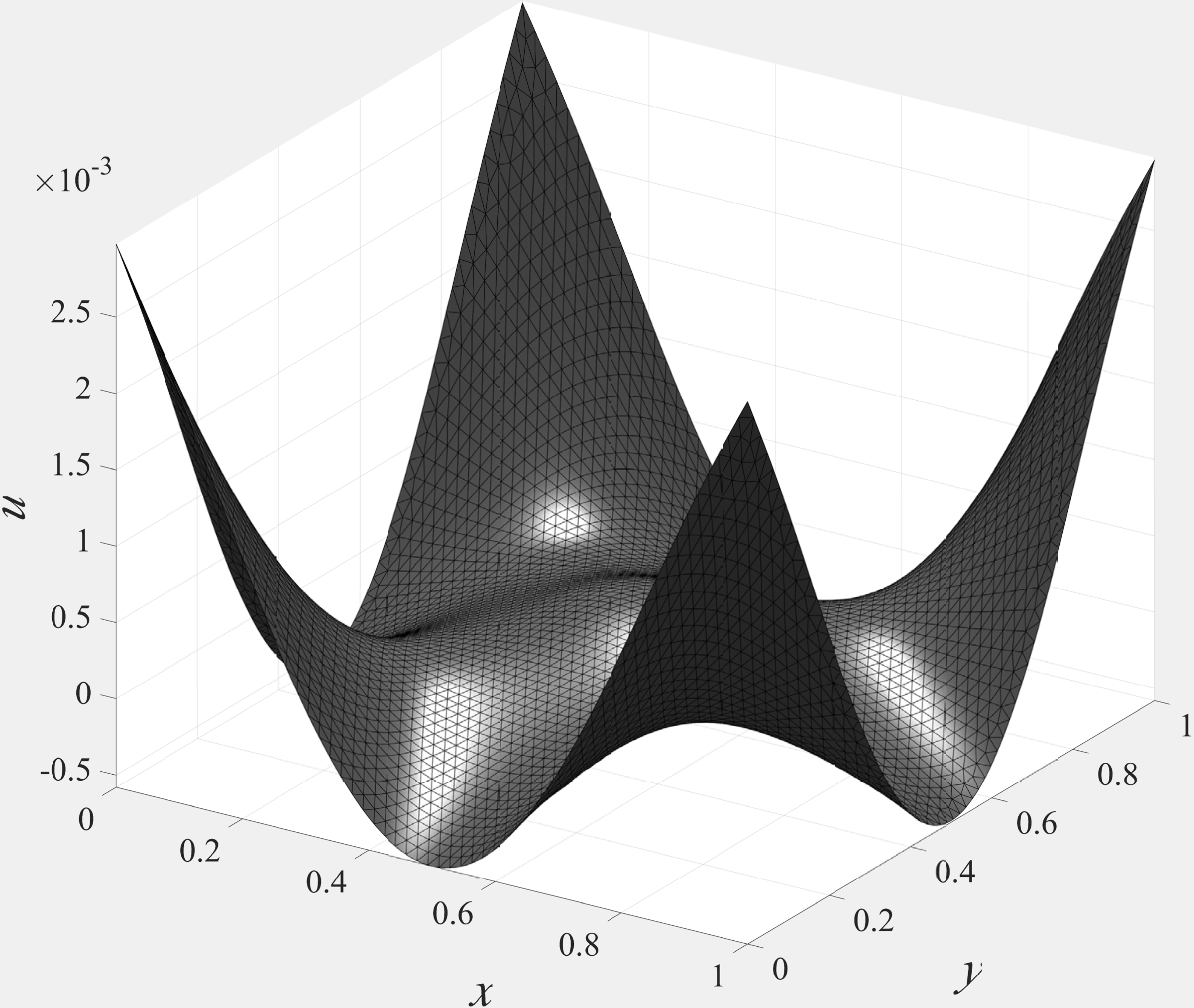}
\includegraphics[scale=0.29]{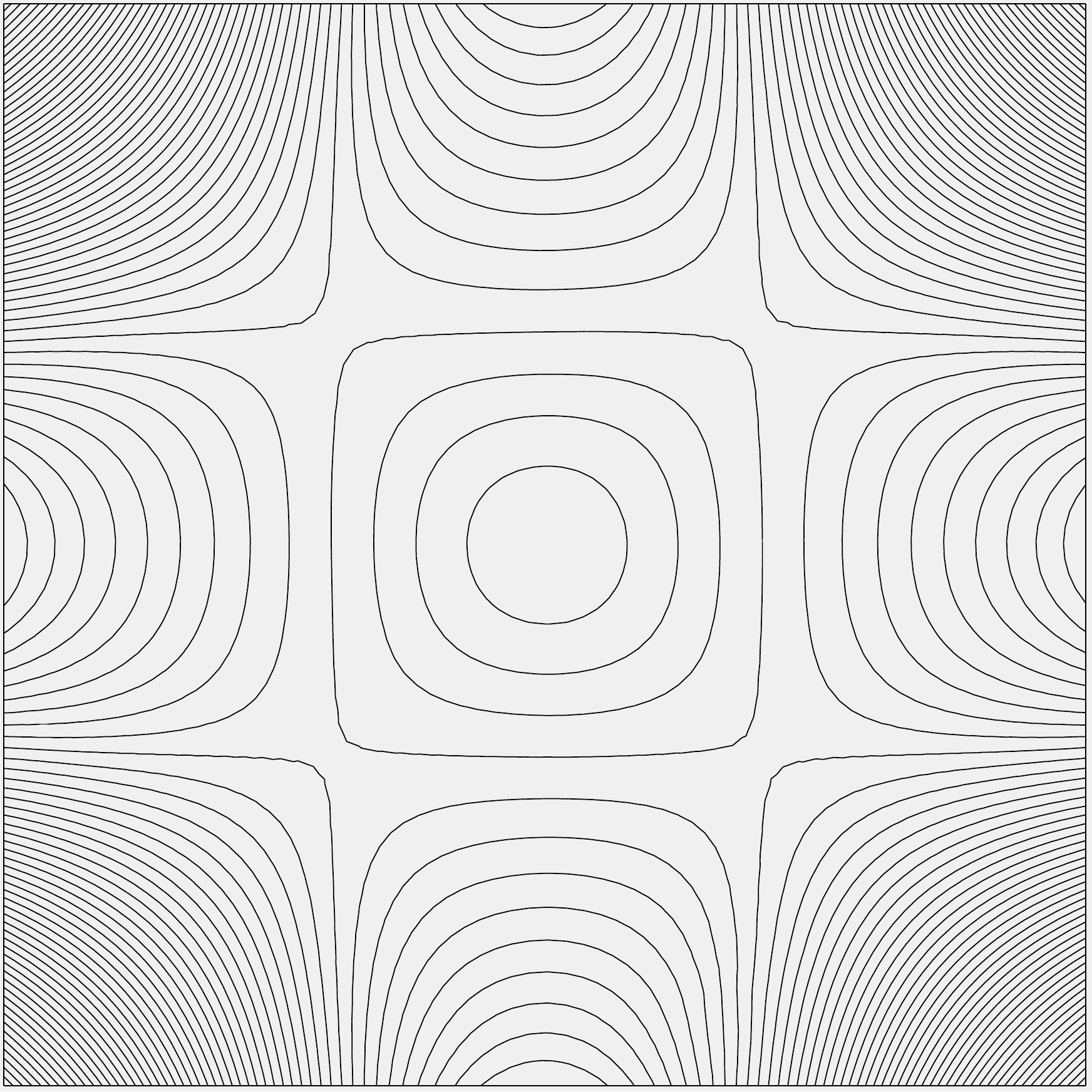}
\end{center}
\caption{Deformations when all beams are clamped.\label{fig:clamped}}
\end{figure}

\begin{figure}
\begin{center}
\includegraphics[scale=0.065]{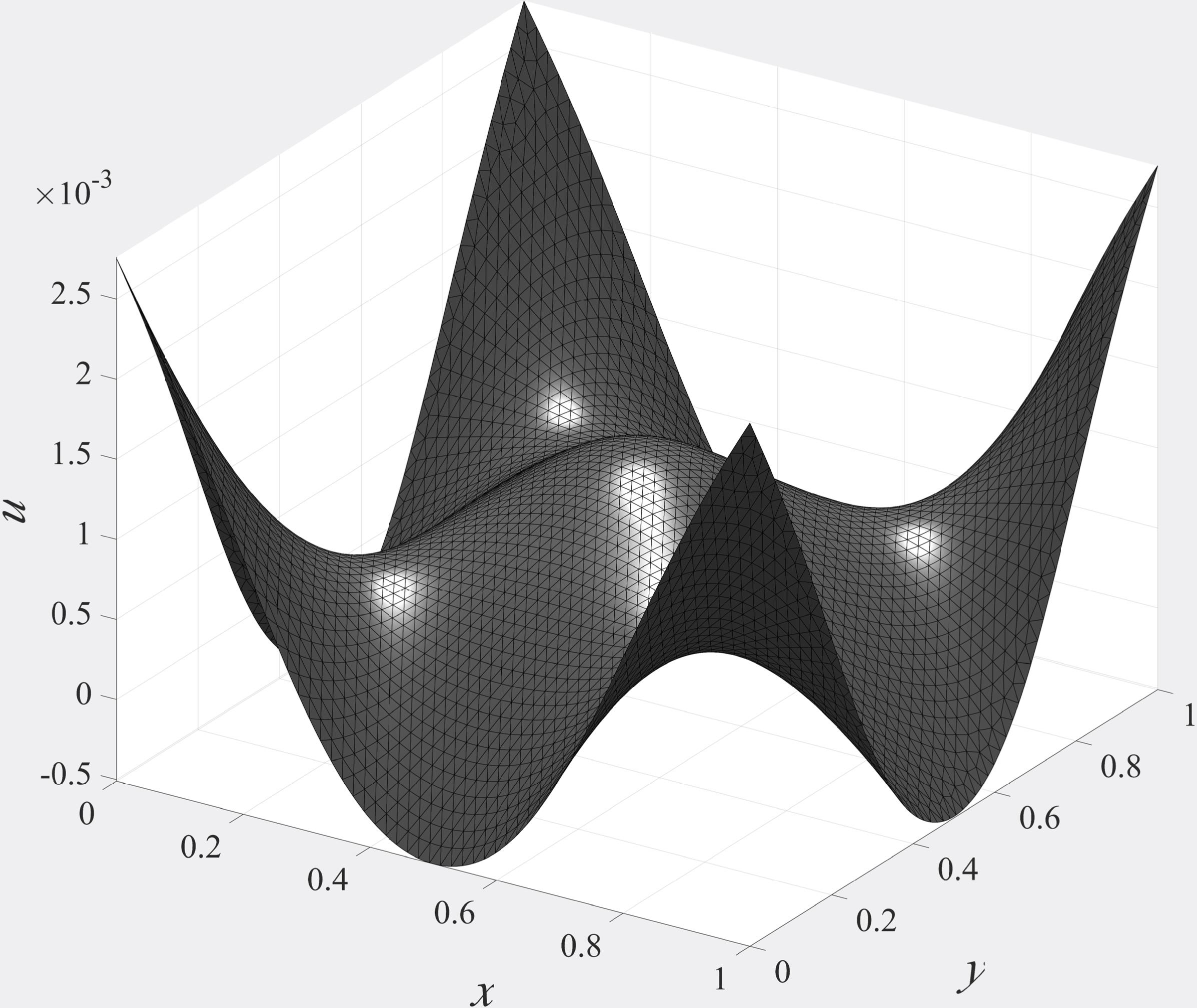}
\includegraphics[scale=0.19]{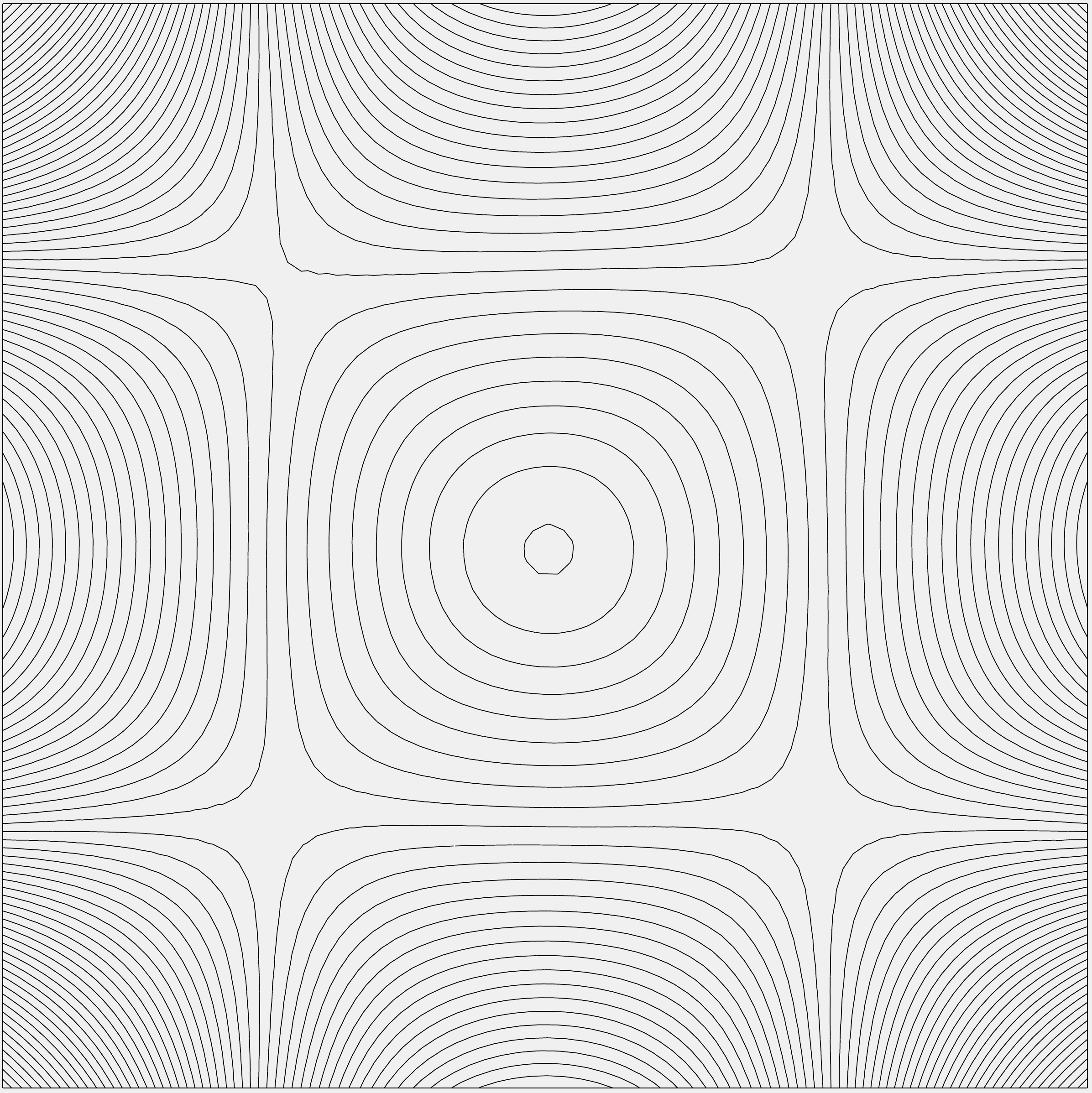}
\end{center}
\caption{Deformations when all beams are simply supported.\label{fig:simple}}
\end{figure}

We first consider the case when the beams are clamped at $x=1$ and free elsewhere.  In Fig. \ref{fig:onesidele} we see the corresponding 
deformation in elevation and isoline plot. Next we consider the case when all beams are clamped, Fig. 
\ref{fig:clamped}, and simply supported, Fig.\ref{fig:simple}. Note the the slight increase in central displacement for the latter.


\section{Conclusions}

We have formulated a continuous/discontinuous Ga\-lerkin method for beam reinforced thin plates. 
The method has the advantage that we can discretize both the beam and plate problem 
with the same standard finite element spaces of continuous piecewise polynomials defined on triangles 
(or quadrilaterals).

\section*{Acknowledgement} This research was supported in part by 
the Swedish Foundation for Strategic Research Grant No.\ AM13-0029, 
the Swedish Research Council Grant No.\ 2013-4708, and 
the Swedish strategic research programme eS\-SEN\-CE. \textcolor{black}{The first author
was supported in part by EPSRC grant EP/P01576X/1.}


\end{document}